\documentclass[10pt]{siamltex}
\usepackage{rotating, booktabs}

\title {Solving nonlinear systems of equations via spectral residual methods: stepsize selection   and applications}

\author{Enrico Meli\footnotemark[1],
Benedetta Morini\footnotemark[2]\  \footnotemark[5],
Margherita Porcelli\footnotemark[3] \footnotemark[5]  \footnotemark[6],
   Cristina Sgattoni\footnotemark[4]\ \footnotemark[5]
 }

\usepackage{amssymb,color,amsmath,latexsym,eucal,cite}
\usepackage{a4wide,graphicx,psfrag} 
\usepackage[latin1]{inputenc}
\usepackage{cases}

\newtheorem{remark}[theorem]{Remark}
\newtheorem{assumption}[theorem]{Assumption}

\newcommand{\bpr}{{\bf Proof.} \hspace{1.5mm}}
\newcommand{\epr}{\hfill $\Box$ \vspace*{1em}}

\def\IR{\hbox{\rm I\kern-.2em\hbox{\rm R}}}
 \def\argmin{\mathop{\rm argmin }}
 \def\liminf{\mathop{\rm liminf }}

 \def\panddot{{\sc Pand}} 
  \def\name{{\sc Srand }} 
    \def\namedot{{\sc Srand}}

\def\Fsigma{{${\tt  F_{bt}}$}} 
\def\Fincr{{ {${\tt F_{in}}$}}} 
\def\Ffmax{{ {${\tt F_{fe}}$}}} 
\def\Fit{{ {${\tt F_{it}}$}}} 
\def\Ibeta{{ {I_{\beta}}}}

\newcommand{\comment}[1]{}
\newcommand{\req}[1]{(\ref{#1})}
\newcommand{\calB}{{\cal B}}

\newcommand{\sigmax}{\sigma_{\max}}
\newcommand{\lmin}{\lambda_{\min}}
\newcommand{\lmax}{\lambda_{\max}}

\newcommand{\bbuno}{\beta_{k,1}}
\newcommand{\bbdue}{\beta_{k,2}}

\newcommand{\di}{v_i}

\newcommand{\lambdai}{\lambda_i}
\newcommand{\muik}{\mu^i_k}
\newcommand{\muikplus}{\mu^i_{k+1}}

\newcommand{\BBu}{{\small{\rm BB1}} }
\newcommand{\BBd}{{\small{\rm BB2}} }
\newcommand{\BBalt}{{\small{\rm ALT}} }
\newcommand{\abbu}{{\small{\rm ABB01}} }
\newcommand{\abbo}{{\small{\rm ABB08}} }
\newcommand{\abbminu}{{\small{\rm ABBm01}} }
\newcommand{\abbmino}{{\small{\rm ABBm08}} }
\newcommand{\dabbmino}{{\small{\rm DABBm}} }

\newcommand{\eqdef}{\stackrel{\rm def}{=}}

\newcounter{algo}[section]
\renewcommand{\thealgo}{\thesection.\arabic{algo}}
\newcommand{\algo}[3]{\refstepcounter{algo}
\begin{center}\begin{figure}[h]
\framebox[\textwidth]{
\parbox{0.95\textwidth} {\vspace{\topsep}
{\bf Algorithm \thealgo : #2}\label{#1}\\
\vspace*{-\topsep} \mbox{ }\\
{#3} \vspace{\topsep} }}
\end{figure}\end{center}}

\begin{document}
\comment{

\begin{titlepage}

\end{titlepage}
}

\date{\documentdate}  

\maketitle

\footnotetext[1]{Dipartimento  di Ingegneria Industriale, Universit\`a degli Studi di Firenze,
via S. Marta 3,  50134 Firenze, Email: enrico.meli@unifi.it}
\footnotetext[2]{Dipartimento  di Ingegneria Industriale, Universit\`a degli Studi di Firenze,
viale G.B. Morgagni 40,  50134 Firenze,  Italia. Email:
benedetta.morini@unifi.it. }
\footnotetext[3]{Dipartimento  di Matematica, AM$^2$, Universit\`a di Bologna,
Piazza di Porta San Donato 5, 40126 Bologna, Italia.  Email:
margherita.porcelli@unibo.it}
\footnotetext[4]{Dipartimento di Matematica e Informatica ``Ulisse Dini'', Universit\`a degli Studi di Firenze,
viale G.B. Morgagni 67a,  50134 Firenze,  Italia. Email:
cristina.sgattoni@unifi.it}
\footnotetext[5]{Member of the INdAM Research Group GNCS.}
\footnotetext[6]{Institute of Information Science and Technologies "A. Faedo", ISTI--CNR, Via Moruzzi 1 Pisa, Italia.}

\begin{abstract}
Spectral residual methods are derivative-free and low-cost per iteration procedures  for solving nonlinear systems of equations.
They are generally coupled with a nonmonotone  linesearch strategy and compare well with Newton-based methods for large nonlinear 
systems and sequences of nonlinear systems.
The residual vector is used as the  search direction and   choosing the steplength  has a crucial impact on the performance.
In this work we address both theoretically and experimentally the steplength selection and provide results on 
a real application such as a rolling contact problem.
\end{abstract}

\section{Introduction}
This work addresses the solution of the nonlinear system  of equations
\begin{equation}\label{ns}
F(x)=0,
\end{equation}
with $F:\mathbb{R}^n\rightarrow \mathbb{R}^n$ continuously differentiable,  by means of spectral residual methods.
Spectral residual methods were introduced in \cite{Sane} and starting from the proposal in \cite{Dfsane}
consist of iterative procedures for solving (\ref{ns}) without the use of derivative information.
Given the iterate $x_k$, these methods use  the residual vectors $\pm F(x_k)$ in a systematic way and select the step 
$x_{k+1}-x_k$ along either the direction 
$(-\beta_k F(x_k))$ or $(\beta_k F(x_k))$ with $\beta_k$ being a nonzero steplength inspired by the Barzilai and Borwein
method for the unconstrained minimization problem $\min_{x\in \mathbb{R}^n} f(x)$.
Similarly to the Barzilai and Borwein method for unconstrained optimization, 
$\|F\|$ does not decrease monotonically along iterations and its effectiveness 
heavily relies on the steplength $\beta_k$ used.

Spectral residual methods have received a large attention since they are
low-cost per iteration and require a low  memory storage being matrix free, see e.g.
\cite{GS, Psane, Sane, Dfsane, Pand, Gas, bumi, Var}. They belong to the class of Quasi-Newton methods which are particularly attractive when the Jacobian  matrix of $F$ is not available analytically or its computation is not relatively easy. Quasi-Newton methods showed to be effective both in the solution of large nonlinear systems and in the solution of sequences of medium-size nonlinear systems as those arising in applications where sequences are generated by model refinement procedures, see e.g.,  \cite{Dfsane, GS, Gas, Sane, Energy, Var}. 

It is well known that the performance of the Barzilai and Borwein method does not depend on the decrease of the objective function
at each iteration but  relies on the relationship between the steplengths used and the eigenvalues of the  average Hessian matrix of $f$
\cite{Birgin, Flet, TesiR}. Based on such feature, several strategies for steplength selection have been proposed to
enhance the performance of the method, see e.g., \cite{DF, DHSZ, Flet, DADSRT, Zanni, FZanni}.
On the other hand, to our knowledge, an analogous study of the relationship between the  steplengths originated by spectral methods 
and the eigenvalues of the  average Jacobian  matrix of $F$ has not been carried out, and the impact of the choice 
of the steplenghts on the convergence history has not been investigated in details.
The aim of this paper is to analyze the properties of the spectral residual steplengths 
and study how they affect the performance of the methods. This aim is addressed both from a theoretical and experimental point of view.

The main contributions of this work are: the theoretical analysis of the steplengths 
proposed in the literature and of their impact on the norm of $F$ also with respect to the nonmonotone behaviour imposed by  
globalization strategies; the analysis of the performance of spectral methods with various rule for updating the steplengths.  
Rules based on adaptive strategies that suitably combine small and large steplengths result by far more effective than rules based on static choices of $\beta_k$ and, inspired by the steplength rules proposed in the literature for unconstrained minimization problems, 
we propose and extensively test adaptive steplength strategies.
Numerical experience is conducted on sequences of nonlinear systems arising  from
rolling contact models which play a central role in many important applications, such as rolling bearings and wheel-rail interaction \cite{Kalker1,Kalker2}. Solving these models gives rise to  sequences  which consist of a large number of medium-size nonlinear systems 
and represent a relevant  benchmark test set for the purpose of this work.


The paper is organized as follows.
Section 2 introduces spectral residual methods. In Section 3 and 4 we provide a theoretical analysis of the  steplengths including their impact on the behaviour of $\|F_k\|$ and on a standard nonmonotone linesearch.
{\color{black} In Section 5 we introduce the spectral residual method used in our tests and provide a theoretical
investigation. The experimental part is developed in Section 6 where we describe several strategies for selecting the steplength, introduce our test set and   discuss the numerical results obtained.
Some conclusions are presented in Section 7.}

\subsection{Notations} 
The symbol $\|\cdot\|$ denotes the Euclidean norm, $I$ denotes the identity matrix, $J$ denotes the Jacobian matrix
of $F$. Given a symmetric matrix $M$,  $\{\lambda_i(M)\}_{i=1}^n$ denotes the set of eigenvalues of $M$,
$\lambda_{\min}(M)$ and $\lambda_{\max}(M)$ denote the minimum and maximum eigenvalue of $M$ 
respectively, and  $\{v_i\}_{i=1}^n$ denotes a set of associated orthonormal eigenvectors. Given a
sequence of vectors $\{x_k\}$, for any function $f$  we let $f_k=f(x_k)$.


\section{Preliminaries}
In the seminal paper  \cite{BB} Barzilai and Borwein proposed a gradient method for the unconstrained minimization 
\begin{equation}\label{minf}
\min_{x\in \mathbb{R}^n} f(x),
\end{equation}
where  $f:\mathbb{R}^n\rightarrow \mathbb{R}$ is a given differentiable function. 
Given an initial guess $x_0\in \mathbb{R}^n$, the Barzilai-Borwein (BB) iteration is defined by
\begin{equation}\label{BB}
x_{k+1}=  x_k-  \alpha_k \nabla f_k,
\end{equation}
where $\alpha_k$ is a positive steplength inspired by Quasi-Newton methods for unconstrained optimization \cite{DS}.
In Quasi-Newton methods, the step $p_k=x_{k+1}-x_k$ solves the linear system
\begin{eqnarray}
& & B_k p_k= -\nabla f_k,\label{QN} 
\end{eqnarray}
and  $B_k$, $k\ge 1$, satisfies the secant equation, i.e., 
\begin{eqnarray}
& & B_k p_{k-1}=z_{k-1}, \quad p_{k-1} =x_{k}-x_{k-1},\quad z_{k-1}=\nabla f_k-\nabla f_{k-1}.\label{secanteq}
\end{eqnarray}
Letting  $B_k= \alpha^{-1}\, I$ and  imposing condition \eqref{secanteq}, Barzilai and Borwein derived two steplengths
which are  the least-square solutions of the following  problems: 
\begin{eqnarray}
&  \alpha_{k,1}&= \argmin_\alpha \| \alpha^{-1} p_{k-1} - z_{k-1} \|_2^2  =\frac{p_{k-1}^T p_{k-1}}{p_{k-1}^T z_{k-1}},\\
& \alpha_{k,2}&= \argmin_\alpha \|p_{k-1} - \alpha z_{k-1} \|_2 ^2   =    \frac{p_{k-1}^T z_{k-1}}{z_{k-1}^T z_{k-1}}.
\end{eqnarray}
The second least-squares formulation  is obtained from the first by symmetry. The steplength $\alpha_k$ in (\ref{BB}) is 
set to be positive, bounded away from zero and not too large, i.e., $\alpha_k\in [\alpha_{\min}, \, \alpha_{\max}]$ for some 
positive $\alpha_{\min}$,  $\alpha_{\max}$; to this end, one of the two 
scalars $\alpha_{k,1}, \, \alpha_{k,2}$ is used and the thresholds $\alpha_{\min}$, $\alpha_{\max}$ are applied to it, see e.g., \cite{Birgin, Flet, Zanni}.

Choosing $B_k=  \alpha^{-1}\, I$  yields  a low-cost iteration while 
the use of the  steplengths $\alpha_{k,1}$, $\alpha_{k,2}$ yields  a considerable 
improvement in the performance with respect to the classical steepest descent method \cite{BB, Flet}. The BB  method
is commonly employed in the solution of large  unconstrained optimization problems (\ref{minf}) and 
the behaviour of the sequence $\{f(x_k)\}$ is typically nonmonotone, possibly severely nonmonotone,
in both the cases of quadratic and general nonlinear functions $f$   \cite{Flet, GHR, Ray2}.
The performance of the BB method depends on the relationship between the steplength $\alpha_k$ and the eigenvalues of the average Hessian matrix 
$\int_0^1 \nabla^2 f(x_{k-1}+t\, p_{k-1}) \,dt$; hence this approach is also denoted as {\em spectral method} 
and an extensive  investigation on  steplength's selection has been carried on \cite{DF, DHSZ, Flet, DADSRT,Zanni, FZanni}.

The extension of this approach to the solution of nonlinear systems of equations (\ref{ns}) 
was firstly proposed by La Cruz and Raydan in  \cite{Sane}. Here we summarize such a proposal and the issues that were inherited by 
subsequent procedures falling into such framework and  designed for both 
general nonlinear systems  \cite{GS, Psane, Sane, Dfsane, Pand, Gas, Var} and for monotone nonlinear systems \cite{ZhangZhou,Yu,Aw,Moh1,Li,Liu}.
Instead of applying the spectral method to the merit function 
\begin{equation}\label{merit}
f(x)=\|F(x)\|^2,
\end{equation}
the BB approach is specialized to the Newton equation
yielding the so-called {\em spectral residual method}. Thus, let  $p_{-}$ satisfy the linear system
\begin{eqnarray}
& & B_k p_{-}= -F_k, \label{QNN}
\end{eqnarray}
and  let $B_k=  \beta^{-1} I$  satisfy the secant equation
$$
B_k p_{k-1}=y_{k-1}, \quad p_{k-1} =x_{k}-x_{k-1}, \quad  y_{k-1}= F_k- F_{k-1}.   
$$

Reasoning as in BB method, two steplengths are derived:
\begin{eqnarray}
\bbuno &=& \frac{p_{k-1}^T p_{k-1}}{p_{k-1}^T y_{k-1}}, \label{beta1}
\\
\bbdue &=& \frac{p_{k-1}^T y_{k-1}}{y_{k-1}^T y_{k-1}}. \label{beta2}
\end{eqnarray}
These scalars may be positive, negative or even null; moreover $\bbuno$ is not well defined if $p_{k-1}^T y_{k-1}=0$
and $\bbdue$ is not well defined if $y_{k-1}=0$.
In practice, the steplength $\beta_k$ is chosen equal  either   to $\bbuno$ or to $\bbdue$ as long as it results 
to be bounded away from zero and $|\beta_k|$ is  not too large, i.e., $|\beta_k|\in [\beta_{\min},\, \beta_{\max}]$ 
for some positive $\beta_{\min}$, $\beta_{\max}$. 
The step  resulting from \eqref{QNN} turns to be of the form 
$p_{-}= -\beta_kF_k.$
But, once $\beta_k$ is fixed, the $k$th  iteration of the spectral residual method employs  the residual directions $\pm F_k$ in a systematic way and tests both 
the steps  
$$
p_{-}= -\beta_kF_k \quad \mbox{and} \quad p_{+}= +\beta_kF_k, 
$$
for acceptance using a suitable linesearch strategy.
The use of both directions  $\pm F_k$ is motivated by the fact that, contrary to $(-\alpha_k \nabla f_k)$, $\alpha_k>0$,  in \req{BB},  $ (-\beta_k F_k)$ is not necessarily a descent direction for
(\ref{merit}) at $x_k$; the value $\nabla f_k^T  (-\beta_k F_k) =-2\beta_k F_k^T J_k F_k$ could be positive, negative or  null. 
On the other hand, if $F_k^T J_k F_k\neq 0$, trivially  either  $ (- \beta_k F_k)$ or  $  \beta_k F_k$ is a descent direction for $f$.

Analogously to the spectral method, the spectral  residual method is characterized by a nonmonotone behaviour of $\{\|F_k\|\}$
and is implemented using nonmonotone line search strategies.
The  adaptation of the spectral method to nonlinear systems is low-cost per iteration since 
the computation of $\bbuno$ and $\bbdue$ is inexpensive
and the memory storage is low, and  turned out to be effective in the solution of medium and large 
nonlinear systems, see e.g.,  \cite{GS, Psane, Sane, Dfsane, Pand,  Var}.

Unlike the context of  BB method for unconstrained optimization, to our knowledge a systematic analysis of the stepsizes $\bbuno$ and $\bbdue$ in the context of the solution of nonlinear systems and their impact on convergence history  has not been carried out. The steplength $\bbuno $ 
has been used  in most of the works on this subject \cite{Psane, Sane, Dfsane, Pand, Gas}. On the other hand,  in \cite{GS} it was 
observed experimentally that alternating $\bbuno$ and $\bbdue$ along iterations was beneficial for the performance and  in \cite{Var} it was 
observed experimentally that using  $\bbdue$ performed better in terms of robustness with respect to using  $\bbuno$.

In the next two sections  we will analyze the two steplengths $\bbuno$ and $\bbdue$ and provide: their expression in terms of the spectrum of average 
matrices associated to the Jacobian matrix of $F$; their mutual relationship;  their  impact   on the behaviour of
$\|F_k\|$ and on a standard nonmonotone linesearch.

The matrices involved in our analysis are the following. 
Given a square matrix $A$, we let  $A_S= \frac12 (A+A^T) $ be  the symmetric part of $A$,  $G_{k-1}$ be the  average  matrix 
associated to the Jacobian $J$ of $F$ around $x_{k-1}$
\begin{equation}\label{G12}
  G_{k-1} \eqdef \int_0^1 J(x_{k-1}+t\, p_{k-1}) \,dt,
\end{equation}
and  $(G_S)_{k-1}$ be the  average matrix associated to the symmetric part $J_S$ of $J$ around $x_{k-1}$
\begin{equation}\label{GS}
(G_S)_{k-1} \eqdef \int_0^1 J_S(x_{k-1}+t\, p_{k-1}) \,dt.
\end{equation}
Moreover, given a symmetric matrix $M$ and a nonzero vector $p$, we employ the Rayleigh quotient  defined as 
\begin{equation}\label{RQ}
 q(M,p)= \frac{p^TMp}{p^Tp},
\end{equation}
 and the following property \cite[Theorem 8.1-2]{Golub} 
\begin{equation}\label{eigenvalues}
\lambda_{\min}(M) \le q(M, p) \le \lambda_{\max}(M). 
\end{equation}

\section{Analysis of the steplengths $\bbuno$ and $\bbdue$}\label{sec::stepl}
We analyze the stepsizes $\bbuno$ and $\bbdue$ given in \eqref{beta1} and \eqref{beta2} making the following assumptions.
\vskip 1pt \noindent
\begin{assumption}\label{A1}
The scalars $\bbuno$ and $\bbdue$  are well defined and nonzero.
\end{assumption}
\begin{assumption}\label{A2}
Given $x$ and $p$, 
$F$ is continuously differentiable in an open convex set $D \subset \mathbb{R}^n$ containing $x +t p $ with $t\in [0,1]$.
\end{assumption}
\vskip 5pt
We note that Assumption \ref{A1} holds whenever $p_{k-1}^Ty_{k-1} \ne 0$.

In the following lemma we analyze the mutual relationship between  the stepsizes  $\bbuno$ and $\bbdue$ and give
their characterization 
in terms of suitable  Rayleigh quotients for the average  matrices in (\ref{G12}) and (\ref{GS}). 
We use repeatedly the property 
\begin{equation}\label{symm}
p^T A p = p^T A_S p,
\end{equation} which holds for any square  matrices $A,$ $A_S= \frac12 (A+A^T)$, and any vector $p$ of suitable dimension.

\vskip 5pt
\begin{lemma}\label{lem:prop}
Let Assumption \ref{A1} hold and Assumption \ref{A2} hold with {\color{black} $x=x_{k-1}$, $p=p_{k-1} = \pm \beta_{k-1} F_{k-1}$}.
The  steplengths $\bbuno$, $\bbdue$ are such that:
\begin{description}
\item{{\em  P1)}} they  have the same sign and $|\bbdue| \le |\bbuno|$; 
\item{{\em P2)}} either  it holds $\bbuno \le \bbdue < 0$ or $0 < \bbdue \le \bbuno$;
\item{{\em P3)}} they take the form 
\begin{equation}\label{bb1}
    \bbuno=\frac{1}{\displaystyle q\big((G_S)_{k-1}, p_{k-1}\big)}{\color{black}=\frac{1}{\displaystyle q\big((G_S)_{k-1}, F_{k-1}\big)} },
\end{equation}
and 
\begin{equation}\label{bb2}
\bbdue=\frac{\displaystyle q\big((G_S)_{k-1}, p_{k-1}\big)}{\displaystyle q(G_{k-1}^T G_{k-1}, p_{k-1})}
{\color{black}=\frac{\displaystyle q\big((G_S)_{k-1}, F_{k-1}\big)}{\displaystyle q(G_{k-1}^T G_{k-1}, F_{k-1})}},
\end{equation}
with $q(\cdot,\cdot)$ being the Rayleigh quotient in \eqref{RQ}, $G_{k-1}$ and $(G_S)_{k-1}$ being the matrices
in (\ref{G12}) and (\ref{GS}), respectively.
\end{description}
\end{lemma}

\bpr
By \eqref{beta1} and \eqref{beta2}, we can write
\begin{eqnarray} \label{cos}
\bbdue &=& \frac{p_{k-1}^T p_{k-1}}{p_{k-1}^T y_{k-1}}  \frac{(p_{k-1}^T y_{k-1})^2}{(y_{k-1}^T y_{k-1})(p_{k-1}^T p_{k-1})}\nonumber  \\ 
&=&  \bbuno \frac{\|p_{k-1}\|^2 \| y_{k-1} \|^2 cos^2\varphi_{k-1}}{\|p_{k-1}\|^2 \|y_{k-1}\|^2} \nonumber \\ 
& = & \bbuno \cos^2\varphi_{k-1}, 
\end{eqnarray}
where $\varphi_{k-1} $ is the angle between $p_{k-1}$ and $y_{k-1}$, and P1)  follows.

Property  P2)  follows as well since  $\bbdue\ne 0$ by Assumption \ref{A1}.

As for property P3), by the Mean Value Theorem \cite[Lemma 4.1.9]{DS} and \eqref{G12} we have
$$
y_{k-1} = F_k - F_{k-1} =  \int_0^1 J(x_{k-1}+t p_{k-1})p_{k-1} \,dt =G_{k-1}p_{k-1}.
$$
Then using \eqref{symm} and (\ref{RQ}), $\bbuno$  takes the form
$$
\bbuno=
  \frac{p_{k-1}^T p_{k-1}}{p_{k-1}^T G_{k-1}p_{k-1}} {\color{black}=  \frac{p_{k-1}^T p_{k-1}}{p_{k-1}^T (G_S)_{k-1}p_{k-1}} }
= \frac{1}{q\big((G_S)_{k-1}, p_{k-1}\big)},
$$
while  $\bbdue$  takes the form
$$
\bbdue= 
 \frac{p_{k-1}^T G_{k-1}  p_{k-1}} {p_{k-1}^T (G_{k-1}^T  G_{k-1})  p_{k-1}} \frac{p_{k-1}^T p_{k-1}}{p_{k-1}^T p_{k-1}} = \frac{q\big((G_S)_{k-1}, p_{k-1}\big)}{q(G_{k-1}^TG_{k-1}, p_{k-1})}.
$$
{\color{black} The rightmost  equalities in \req{bb1} and \req{bb2} easily follow using the form of the step $p_{k-1} = \pm \beta_{k-1} F_{k-1}$.}
\epr

The above characterization P3) allows to derive  bounds on the stepsizes $\bbuno$ and $\bbdue$ diversifying cases according 
to the spectral properties of the  Jacobian matrix and the average matrices in (\ref{G12}) and (\ref{GS}).  
 The  relationship between $\bbuno$ and the spectral information of the symmetric part of  average  matrix (\ref{G12}) was observed in 
\cite{Dfsane, Sane, Pand} but the following results are not contained in such references.

\vskip 5pt
\begin{lemma}\label{lem:bounds}
Let Assumption \ref{A1} hold and Assumption \ref{A2} hold with $x=x_{k-1}$, $p=p_{k-1}$. Then, the steplengths $\bbuno$ and $\bbdue$ are such that:
\begin{description}
\item{{\em (i)}} If  the Jacobian $J$ is symmetric and positive definite on the line segment in between $x_{k-1}$ and $x_{k-1}+ p_{k-1}$ then  
$\bbuno$ and $\bbdue$ are positive and
\begin{equation}\label{bJpd}
\frac{1}{\lambda_{\max}(G_{k-1})}
\le \bbdue \le \bbuno \le 
\frac{1}{\lambda_{\min}(G_{k-1})};
\end{equation}
\item{{\em(ii)}} if $(G_S)_{k-1}$ in (\ref{GS}) is positive definite, then $\bbuno$ and $\bbdue$ are positive and 
\begin{equation}\label{bGspd}
\max\bigg\{\frac{1}{\lambda_{\max}\big((G_S)_{k-1}\big)}, \, \bbdue\bigg\}\le \bbuno \le \frac{1}{\lambda_{\min}\big((G_S)_{k-1}\big)},
\end{equation}
\begin{equation}\label{bGspd2}
\frac{\lambda_{\min}\big((G_S)_{k-1}\big)}{\lambda_{\max}(G_{k-1}^TG_{k-1})}\le \bbdue \le \min\bigg\{\frac{\lambda_{\max}\big((G_S)_{k-1}\big)}{\lambda_{\min}(G_{k-1}^TG_{k-1})}, \bbuno\bigg\};
\end{equation}
\item{{\em(iii)}}  if  $(G_S)_{k-1}$ in (\ref{GS}) is indefinite  and  $G_{k-1}$ in \eqref{G12} is nonsingular, then
\begin{description}
\item{{\em (iii.1)}} $\bbuno$ satisfies either
\begin{equation}\label{bb1indef}
 \bbuno \le \min\left\{  \frac{1}{\lambda_{\min}\left((G_S)_{k-1}\right)}, \bbdue \right\} 
\ \ \:  \mbox{or} \ \ \: \bbuno\ge \max \left\{ \frac{1}{\lambda_{\max}\left((G_S)_{k-1}\right)}, \bbdue \right\} ;
\end{equation}
\item{{\em(iii.2)}} $\bbdue$ satisfies either
\begin{equation}\label{bb2indef}0 < \bbdue \le \min\bigg\{\frac{\lambda_{\max}\big((G_S)_{k-1}\big)}{\lambda_{\min}(G_{k-1}^TG_{k-1})}, \bbuno\bigg\},
\end{equation}
or
\begin{equation}\label{bb2indef2}\max\bigg\{\frac{\lambda_{\min}\big((G_S)_{k-1}, \big)}{\lambda_{\max}(G_{k-1}^TG_{k-1})},\bbuno\bigg\}\le \bbdue < 0.
\end{equation}
\end{description}
\end{description}
\end{lemma}

\bpr Consider properties P1), P2) and P3) from Lemma \ref{lem:prop}.
\begin{description}
\item{{ (i)}}
Steplengths $\bbuno$ and $\bbdue$ are positive due to \eqref{bb1}, \eqref{bb2}. The rightmost inequality of (\ref{bJpd}) follows from
(\ref{bb1}) and \eqref{eigenvalues}. The remaining part of (\ref{bJpd}) is proved observing that \eqref{bb2} yields
\begin{equation}\label{bb2_2}
\bbdue = \frac{p_{k-1}^T G_{k-1}^{1/2} G_{k-1}^{1/2} p_{k-1}}{p_{k-1}^T G_{k-1}^{1/2} G_{k-1} G_{k-1}^{1/2} p_{k-1}}= \frac{1}{q(G_{k-1}, G_{k-1}^{1/2} p_{k-1})},
\end{equation}
and using P2) and \eqref{eigenvalues}.
\item{{ (ii)}}
Using \eqref{bb1},\eqref{eigenvalues} and P2) we get positivity of  $\bbuno$ and \eqref{bGspd}.
Consequently, $\bbdue$ is  positive by property  P1), and 
bounds \eqref{bGspd2} can be derived using \eqref{bb2}, \eqref{eigenvalues} and  item P2) of Lemma \ref{lem:prop}.
\item{{ (iii)}}
If $(G_S)_{k-1}$ is indefinite then  its extreme eigenvalues have opposite sign, i.e.,  $\lambda_{\min}\big((G_S)_{k-1}\big) < 0$
and $\lambda_{\max}\big((G_S)_{k-1}\big) >0$.  Hence, \eqref{bb1}, \eqref{eigenvalues} and P2) give \eqref{bb1indef}.
Moreover, since ${G}_{k-1}^T G_{k-1}$ is symmetric and positive definite, we can use, as before, P1) and  \eqref{eigenvalues} and get
\eqref{bb2indef} and \eqref{bb2indef2}.

\end{description}
\epr
\begin{remark}
Lemma \ref {lem:bounds} easily extends to the case where  matrices   are  negative definite.

Item  {\em (ii)} of Lemma \ref{lem:bounds} includes the case where $F$ is strictly monotone, i.e., $(F(x)-F(y))^T(x-y)> 0$ for any $x, y\in \mathbb{R}^n$ with $x\neq y$,
see  e.g.  \cite{Fac}. 
\end{remark}

\section{On the impact of  the steplength $\beta_k$ on  $ \|F_{k+1}\| $}\label{sec::stepF}

In this section we investigate   how the choice of the steplength $\beta_k$ may affect $\|F_{k+1}\|$ in a spectral residual method.
Results are first derived using a generic $\beta_k$  and  discussed  thereafter with respect to the choice of either $\bbuno$ 
 or  $\bbdue$.

The first result concerns  the case where  $J$ is symmetric and analyzes  the residual vector $F_{k+1}$ componentwise.
It heavily relies on the existence of a set of orthonormal eigenvectors for the average matrix $G_k$.   
\vskip 5pt
\begin{lemma}\label{prop::sym}
Suppose that Assumption \ref{A2} holds with $x=x_{k}$ and $p=p_{k}$ and that  the Jacobian $J$ is symmetric. 
Let $p_k=  p_{-} = - \beta_k F_k\neq 0$, $x_{k+1} =x_k+p_k$,  
$\big\{\lambdai\big(G_k\big)\big\}_{i=1}^n$ be the eigenvalues of  matrix $G_k$ in (\ref{G12}) 
and $\{\di\}_{i=1}^n$ be a set
of associated orthonormal eigenvectors. Let  $F_k$ and $F_{k+1}$  be expressed as
$$
F_k  =  \sum_{i=1}^n \muik \di, \qquad F_{k+1} = \sum_{i=1}^n \muikplus \di,\nonumber\\
$$
where $\muik, \muikplus$, $i=1,\dots,n$, are scalars. Then
\begin{eqnarray}
& & F_{k+1} =  (I  - \beta_k G_k) F_k,\label{Fkplus}\\ 
& & \muikplus = \mu^i_{k}\big(1- \beta_k \lambdai(G_k)\big) , \qquad i=1,\dots,n.\label{mkplus}
\end{eqnarray}
Moreover, it  holds:
\begin{description}
\item{\em {(a)}}  if $\beta_k \lambdai (G_k) = 1$, then $|\muikplus| =0$;
\item{\em {(b)}} if $0 < \beta_k \lambdai (G_k) < 2$, then $|\muikplus| < |\muik|$; otherwise  $|\muikplus| \ge  |\muik|$.
\end{description}
\end{lemma}
\vskip 5pt
\bpr
The Mean Value Theorem \cite[Lemma 4.1.9]{DS} gives 
\begin{eqnarray*}
F_{k+1} &=&F_k+ \int_0^1 J(x_k+t p_k) p_k \,dt ,\\
				\end{eqnarray*}
and $p_k=-\beta_k F_k$ and (\ref{G12}) yield  \eqref{Fkplus}. 
Moreover, since $\{\di\}_{i=1}^n$ are orthonormal we have for $i=1,\dots,n$
\begin{eqnarray*}
\muikplus &=& (\di)^T F_{k+1} \\ 
        &=& (\di)^T (I  - \beta_k G_k ) F_k \\
        &=& \muik \big(1-\beta_k \lambdai(G_k)\big),
\end{eqnarray*}
i.e., equation \eqref{mkplus}. 
Consequently, Item (a)  follows trivially; Item   (b)  follows noting that 
$\big|1- \beta_k \lambdai (G_k)\big|<1$ if and only if  $0 < \beta_k \lambdai(G_k) < 2$.
\epr
\begin{remark}
Lemma \ref {prop::sym} trivially extends to the case where  $p_k=p_{+}=\beta_k F_k$.
\end{remark}
\vskip 5pt
If the nonlinear system \eqref{ns} represents the first-order optimality condition of the optimization problem \eqref{minf}
where $f(x)=\frac{1}{2}x^TAx - b^Tx$ is quadratic and $A$ is symmetric and positive definite, then 
the previous lemma reduces to well known results on  the behaviour of the gradient method in terms of the spectrum of the Hessian matrix $A$,
see  \cite{TesiR}.  In fact, the nonlinear residual is $F(x) = Ax-b$ and its Jacobian is constant $J(x)=A, \, \forall x$. Then the following  strict relationship between $F_k$ and the $i$th eigenvalue $\lambda_i(A)$ of the Jacobian holds throughout the iterations
$$
\muikplus = \muik (1-\beta_k\lambda_i(A))=\mu_0^{i} \prod_{j=0}^k(1-\beta_j\lambda_i(A)),
$$
where $\muikplus$ and $\muik$,  $i=1, \ldots n$, are the eigencomponents of 
$F_{k+1}$ and $F_k$ respectively,
with respect to the eigendecomposition of $A$.
As a consequence, a small steplength $\beta_k$, i.e., close to $1/\lambda_{\max}(A)$, can significantly reduce 
the values {\color{black}$|\muikplus|$} corresponding to large eigenvalues $\lambda_{i}(A)$ while 
a small reduction is expected for the scalars {\color{black}$|\muikplus|$} corresponding to small eigenvalues $\lambda_{i}(A)$.
On the contrary, a large steplength $\beta_k$, i.e., close to $1/\lambda_{\min}(A)$, can significantly reduce 
the values {\color{black}$|\muikplus|$} corresponding to small eigenvalues $\lambda_{i}(A)$ while tends to increase  the scalar {\color{black}$|\muikplus|$}
corresponding to large eigenvalues $\lambdai(A)$.
This offers some intuition for  choosing the steplengths by 
alternating  in a balanced way small and large steplengths in order to reduce the eigencomponents, see e.g., \cite[p. 178]{Zanni}.

On the other hand, if $F$ is a general nonlinear mapping then  $G_k$  changes at each iteration and Lemma \ref{prop::sym} 
suggests the expected change of $ F $  from iteration $k$ to iteration  $k+1$ and the following guidelines.
The first guideline concerns the case where $J$ is positive definite. 
A nonmonotone behaviour  of the sequence $\{\|F_k\|\}$ is expected. By Item   (i)  of Lemma \ref{lem:bounds}, both 
$\bbuno$ or $\bbdue$  are  positive and
 $\beta_k \lambdai(G_k)$  lies in the interval 
$\displaystyle \left[\frac{\lambda_{i}(G_{k})}{\lambda_{\max}(G_{k-1})}, \,\frac{\lambda_{i}(G_{k})}{\lambda_{\min}(G_{k-1})}  \right]$ for $i=1,\dots,n$. 
Assuming without loss of generality that the eigenvalues are numbered in nondecreasing order, 
by standard arguments on perturbation theory for the eigenvalues it holds
$$
|\lambdai(G_{k})-\lambdai(G_{k-1})| \le \|G_k-G_{k-1}\|, 
$$
$i=1, \ldots, n$, \cite[Theorem 8.1-6]{Golub}. Thus, if  the Jacobian is Lipschitz continuous 
in an open convex set containing $x_{k-1}+t p_{k-1}$ and $x_k+t p_k$
with constant $L_J>0$, it follows 
$$
\|G_k-G_{k-1}\|\le \frac{L_J}{2}\bigg(\|p_{k-1}\|+\|p_k\|\bigg).
$$
Hence, if $\|p_{k-1}\|$ and/or $\|p_k\|$ are large, by Item (b) no decrease of $\muikplus$ may occur. On the contrary,
for small values of $\|p_{k-1}\|$ and $\|p_k\|$, as occurs if $\{x_k\}$ is convergent, $G_k$ undergoes small changes with respect to 
$G_{k-1}$ and the behaviour of $\muikplus$ shows similarities with the case where $J$ is constant. 
Thus, a small steplength $\beta_k$ close to $1/\lambda_{\max}(G_{k-1})$ can significantly reduce 
the scalars {\color{black}$|\muikplus|$} corresponding to large eigenvalues $\lambda_{i}(G_k)$, while 
a small reduction is expected for the values {\color{black}$|\muikplus|$} corresponding to small eigenvalues $\lambda_{i}(G_k)$.
A large steplength $\beta_k$  close to $1/\lambda_{\min}(G_{k-1})$ can significantly reduce 
the scalars {\color{black}$|\muikplus|$} corresponding to small eigenvalues $\lambda_{i}(G_k)$ 
while tends to increase  the eigencomponents {\color{black}$|\muikplus|$} corresponding to large eigenvalues $\lambdai(G_k)$.
As for the case of a constant Jacobian, these features suggest to choose the steplengths by 
alternating  in a balanced way small and large steplengths in order to reduce the eigencomponents.


 The second guideline concerns the case where $J$ is indefinite and  $\lambda_{\min}(G_k) < 0 < \lambda_{\max}(G_k )$. If 
$\beta_k >0$,  from Item (b)  it follows that  $|\muikplus|$ corresponding to positive $\lambda_{i}(G_k)$ are smaller than 
$|\muik|$ if  $\beta_k \lambdai(G_k)$ is small enough while all  $|\muikplus|$ corresponding to negative eigenvalues 
increase  with respect to  $|\muik|$ and the amplification depends on the magnitude of $\beta_k \lambdai(G_k)$.
If $\beta_k <0$ similar conclusions hold.
In general, a nonmonotone behaviour  of the sequence $\{\|F_k\|\}$ is expected
but a possibly large increase of $\|F_{k+1}\|$ with respect to $\|F_k\|$ does not occur if 
$\{|\beta_k \lambdai(G_k)|\}_{i=1, \ldots,n}$
are small or of moderate size.
Since a small value of  $\{|\beta_k \lambdai(G_k)|\}_{i=1, \ldots,n}$ might 
be induced by a small value of $|\beta_k|$, the use of $\bbdue$ might be advisable taking into account that $|\bbdue|\le |\bbuno|$ and
$\bbuno$ can arbitrarily grow in the indefinite case (see Lemma \ref{lem:bounds}).

\subsection{On the impact of the steplength $\beta_k$ in the approximate norm descent  linesearch}\label{sec::stepnm}
In this section we embed the spectral residual  method in a general globalization scheme based on the  so-called approximate norm descent condition \cite{LiFuk}
\begin{equation}\label{normdes}
\| F_{k+1}\| \le (1+\eta_k) \|F_k\|,
\end{equation}
where $\{\eta_k\}$ is a positive sequence  satisfying
\begin{equation}\label{etaserie}
\sum_{k=0}^{\infty} \eta_k < \eta < \infty.
\end{equation}
Intuitively, large values of $\eta_k$ allow a highly nonmonotone behaviour of
$\|F_k\|$ while small values of $\eta_k$ promote  the decrease of $\|F \|$.
Several linesearch strategies in the literature fall in this scheme
 \cite{LiFuk,Pand,Gas,Gonc}. The main idea is that, given $x_k$, the steps take the form
\begin{equation}\label{iter2}
 p_- = -\gamma_k \beta_k F_k  \quad \text{or} \quad p_+ = +\gamma_k \beta_k F_k
\end{equation} 
where the sign $\pm$ and $\gamma_k \in (0,1]$ are selected so that  \req{normdes} is satisfied.
The scalar $\gamma_k$ can be computed using  a
backtracking process.
Enforcing condition (\ref{normdes})  ensures the convergence of the sequence $\{\|F_k\|\}$ \cite[Lemma 2.4]{LiFuk}. 

We now analyse the properties of $\|F_{k+1}\|$ as a function of the stepsize $\gamma_k\beta_k$ 
and determine conditions on $\gamma_k\beta_k$ which enforce  \eqref{normdes}. 
First of all we observe that by the Mean Value Theorem \cite[Lemma 4.1.9]{DS} and (\ref{iter2}) we have  
\begin{equation}\label{Fkplusgen}
F_{k+1} = (I\pm\gamma_k\beta_kG_k)F_k.
\end{equation}
Using this equation we can write
\begin{equation}
\|F_{k+1}\|^2=\|F_k\|^2\pm2\gamma_k\beta_k F_k^T (G_S)_k F_k+\gamma_k^2\beta_k^2F_k^TG_k^TG_kF_k, \label{Fplgen}
\end{equation}
and analyze the fulfillment of either the decrease of $\|F\|$ or (\ref{normdes}) as given below. 
\vskip 5pt
\begin{theorem}\label{thm:unsym}
Suppose that Assumption \ref{A1} holds and Assumption \ref{A2} holds with $x=x_k$ and $p=p_k$. 
Suppose $F_k^T J_kF_k\neq 0$  
and $F_k^TG_kF_k\neq 0$ 
 with $G_k$ given in (\ref{G12}). Let $\Delta=q\big((G_S)_k, F_k\big)^2+(\eta_k^2+2 \eta_k)q(G_k^T G_k, F_k)$, then
\begin{description}
\item{(1)} If $x_{k+1} =x_k+p_k$, $p_k= p_-=- \gamma_k\beta_k F_k$,  $\gamma_k\in (0,1]$, we have that $\|F_{k+1}\| < \|F_k\|$
when 
\begin{eqnarray}
& \beta_kq\big((G_S)_k, F_k\big)>0 & \ \mbox{ and } \
\gamma_k \big|\beta_k\big| < 2\,  \frac{\big| q\big((G_S)_k, F_k\big)\big|}{q(G_k^T G_k, F_k)} . \label{dlin1}
\end{eqnarray}
Condition  (\ref{normdes})
is satisfied when 
\begin{eqnarray}
& &  \frac{  q\big((G_S)_k, F_k\big) -\sqrt{ \Delta}}{q(G_k^T G_k, F_k)} \le \gamma_k  \beta_k \le  \frac{  q\big((G_S)_k, F_k\big)  +\sqrt{ \Delta}}{q(G_k^T G_k, F_k)}.  \label{dlin3}
\end{eqnarray}
\item{(2)} If $x_{k+1} =x_k+p_k$, $p_k= p_+ = \gamma_k\beta_k F_k$,  $\gamma_k\in (0,1]$, we have that $\|F_{k+1}\| < \|F_k\|$
when 
\begin{eqnarray}
& \beta_kq\big((G_S)_k, F_k\big)<0 & \ \mbox{ and } \
\gamma_k \big|\beta_k\big| < 2\,  \frac{\big| q\big((G_S)_k, F_k\big)\big|}{q(G_k^T G_k, F_k)} 
\label{dlin2}
\end{eqnarray}
Condition  (\ref{normdes})
is satisfied when 
\begin{eqnarray}
& &  \frac{  -q\big((G_S)_k, F_k\big) -\sqrt{ \Delta}}{q(G_k^T G_k, F_k)} \le \gamma_k  \beta_k \le  \frac{  -q\big((G_S)_k, F_k\big)  +\sqrt{ \Delta}}{q(G_k^T G_k, F_k)}.  \label{dlin4}
\end{eqnarray}
\end{description}
\end{theorem}
\bpr 
Concerning Item (1), using  (\ref{Fkplusgen}) we get
\begin{eqnarray*}
\|F_{k+1}\|^2 
            & = &    {\color{black}  \|(I -\gamma_k\beta_kG_k)F_k \|^2}\\
                  & = & \Big(1 -2 \gamma_k\beta_k \frac{F_{k}^T  (G_S)_k F_{k}}{\|F_k\|^2}   + \gamma_k^2 \beta_k^2 \frac{ F_{k}^T  G_k^T G_k F_{k}}{\|F_k\|^2}\Big )\|F_k\|^2\nonumber \\
                  & = & \Big(1-2\gamma_k \beta_k q\big((G_S)_k, F_k\big)+\gamma_k^2 \beta_k^2 q(G_k^T G_k, F_k)\Big)\|F_k\|^2.    
\end{eqnarray*}
Noting that by assumption $q\big((G_S)_k, F_k\big)\neq 0$ and $q(G_k^T G_k, F_k)>0$,  $\|F_{k+1}\|<\|F_k\|$ holds if 
$$
\beta_k q\big((G_S)_k, F_k\big)>0  \quad \mbox{ and} \quad -2\gamma_k \beta_k q\big((G_S)_k, F_k\big)+\gamma_k^2 \beta_k^2 q(G_k^T G_k, F_k)<0,
$$
and these conditions can be rewritten as in (\ref{dlin1}). Condition (\ref{dlin3}) follows trivially.

Item $(2)$ follows analogously. 
From (\ref{Fkplusgen}) and imposing 
and  $\|F_{k+1}\|<\|F_k\|$ we get the condition  
$$
\beta_k q\big((G_S)_k, F_k\big)<0  \quad \mbox{ and} \quad  2\gamma_k \beta_k q\big((G_S)_k, F_k\big)+\gamma_k^2 \beta_k^2 q(G_k^T G_k, F_k)<0
$$
which is equivalent to (\ref{dlin2}). Condition (\ref{dlin4}) follows trivially.
\epr 
\vskip 5pt 
We remark that, due to the form of $G_k$ and $(G_S)_k$, 
conditions (\ref{dlin1})--(\ref{dlin4}) are implicit in $\gamma_k\beta_k$.
The above theorem supports testing the   two steps (\ref{iter2}) systematically because of the following fact.
At $k$-th iteration, $\beta_k$, $q\big( J_k, F_k\big)$ and $q(J_k^T J_k, F_k)$ are given and 
 by continuity of the Jacobian, the Rayleigh quotients   $q\big((G_S)_k, F_k\big)$ and 
$ q(G_k^T G_k, F_k)$ tend to $q\big(J_k, F_k\big)$ and $q(J_k^T J_k, F_k)$ respectively 
as  $\gamma_k$ tends to zero. 
Hence,  if $\gamma_k$ is sufficiently small then  
$$ 
\frac{q\big(J_k, F_k\big)-\epsilon}{q\big(J_k^T J_k, F_k\big)+\epsilon}
 \le \frac{q\big((G_S)_k, F_k\big)}{q\big(G_k^T G_k, F_k\big)} \le  \, \frac{q\big(J_k, F_k\big)+\epsilon}{q\big(J_k^T J_k, F_k\big)-\epsilon},
$$
and if  $ 0 <\epsilon<\frac 1 2 \min\{\left |q\big( J_k, F_k\big) \right |,\,q(J_k^T J_k, F_k) \}$ then $\frac{q\big((G_S)_k, F_k\big)}{q\big(G_k^T G_k, F_k\big)}$ has the same sign as 
$ \frac{q\big(J_k, F_k\big)}{{q\big(J_k^T J_k, F_k\big)}}$. 
Consequently, for $\gamma_k$ sufficiently small, either  condition  (\ref{dlin1}) or (\ref{dlin2}) is fulfilled.
Analogous considerations can be made for conditions (\ref{dlin3}) and (\ref{dlin4}).

As a final comment, the previous theorem suggests that a small {\color{black}$|\beta_k|$  promotes the fulfillment of conditions} (\ref{dlin1}) and (\ref{dlin2}) or (\ref{dlin3}) and (\ref{dlin4}). 
Again, by   Lemma \ref{lem:bounds}, the use of $\bbdue$ may be advisable taking into account that $|\bbdue|\le  |\bbuno|$ and that 
$\bbuno$ can arbitrarily grow in the indefinite case; taking the steplength equal to $\bbuno$  may cause a large number of backtracks
and an erratic behaviour of $\{\|F_k\|\}$ as long as $\eta_k$ is sufficiently large.

%
%
{\color{black}
\section{A spectral residual approximate norm descent method}\label{sec:srand}
In this section we describe a spectral residual algorithm which implements a line-search along $\pm F_k$ 
and enforces the approximate norm descent condition \req{normdes}.
We also discuss the convergence properties of the method and provide sufficient conditions for the convergence of the sequence
$\{\|F_k\|\} $ to zero.

The Projected Approximate Norm Descent (\panddot) algorithm  was developed in \cite{Pand} 
for solving convexly constrained nonlinear systems. Among its  variants  proposed in \cite{Pand, Gas} and based on Quasi-Newton methods, 
we consider the spectral residual implementation for unconstrained nonlinear systems which is the focus of this work
and denote it as Spectral Residual Approximate Norm Descent (\namedot) method.

Given the current iterate $x_k$, a  new iterate $x_{k+1}$ is computed as
$x_{k+1}= x_k + p_k$ with $p_k$ given by either $(- \gamma_k \beta_k F_k)$ or
$(+ \gamma_k\beta_k F_k)$, $\gamma_k \in (0,1]$.
The main phases of \name are as follows. First, the scalar $\beta_k$ is chosen to that  $|\beta_k|\in [\beta_{\min}, \beta_{\max}]$.
Second,  the scalar $\gamma_k \in (0,1]$ is fixed using a backtracking strategy so that either the  linesearch condition
\begin{equation}\label{lin1}
\|F(x_k + p_k)\| \le \big(1-\rho(1+\gamma_k)\big)\|F_k\|,
\end{equation}
holds or the  linesearch condition 
\begin{equation}\label{lin2}
\|F(x_k  + p_k) \| \le (1+\eta_k-\rho \gamma_k) \|F_k\|,
\end{equation}
holds where  $\rho \in (0,1)$ is quite small \cite{DS,Pand} and   $\{\eta_k\}$  is  a positive sequence satisfying   \req{etaserie}.
The linesearch conditions  \eqref{lin1} and \eqref{lin2} are derivative-free; the first condition
imposes at each iteration a sufficient decrease in $\|F\|$ 
which can be accomplished for suitable values of $\pm \gamma_k\beta_k F_k$ as long as  $F_k^T J_kF_k\neq 0$, and is crucial 
for  establishing  results on the convergence of $\{\|F_k\|\}$ to zero.
On the other hand, the second condition  allows for an increase of $\|F\|$ depending on the magnitude of $\eta_k$.
Trivially, \eqref{lin1} implies \eqref{lin2} and both imply the approximate norm descent condition \req{normdes}.

The formal description of the \name method is reported in  Algorithm \ref{pand_algo}
where we deliberately do not specify the form of the stepsize $\beta_k$. Termination of Step 2 is guaranteed by Theorem \ref{thm:unsym}.
The theoretical properties of \name given in \cite[Theorem 4.2 and Theorem 4.3]{Pand} are summarized in the following theorem.

\begin{theorem}\label{convF}  Let the positive sequence $\{\eta_k\}$ satisfy (\ref{etaserie}) and 
let $\{x_k\}$ be the sequence generated by the \name algorithm. Then
\begin{enumerate}
\item the sequence $\{x_k\}$ is convergent and consequently the sequence $\{\|F_k\|\}$ is convergent;
\item  the sequence $\{\gamma_k \|F_k\|\}$ is convergent and such that 
$\lim_{k\rightarrow \infty} \gamma_k\|F_k\|=0$;
\item if (\ref{lin1}) is satisfied for infinitely many $k$, then $\lim_{k\rightarrow \infty}\|F_k\|=0 $. 
\end{enumerate}
\end{theorem}
\algo{pand_algo}{The \name algorithm}{
Given   $x_0\in \mathbb{R}^{n}$, $0 < \beta_{\min}< \beta_{\max}$, $\beta_0 \in  [\beta_{\min}, \beta_{\max}]$, $\rho, \, \sigma \in (0,1)$, 
a positive sequence $\{\eta_k\}$ satisfying (\ref{etaserie}).
\vskip 4pt
\noindent
If $\|F_0\| =0$ stop.  \\ 
For $k=0,\,1, \,2,\, \ldots\,\,$  do\\
\hspace*{15pt} 1. Set $\gamma=1$.\\
\hspace*{15pt} 2. Repeat\\
\hspace*{44pt} 2.1~ Set  $p_{-}= - \gamma \beta_k F_k$ and $p_{+}= \gamma \beta_k F_k$.\\
\hspace*{44pt} 2.2~ If $p_{-}$  satisfies (\ref{lin1}), set $p_k=p_{-}$ and  go to Step 3.\\
\hspace*{44pt} 2.3~  If  $p_{+}$ 
                    satisfies (\ref{lin1}), set $p_k=p_{+}$ and  go to Step 3.\\
\hspace*{44pt} 2.4~ If 
                 $p_{-}$  satisfies  (\ref{lin2}), set  $p_k=p_{-}$ and go to Step 3.\\
\hspace*{44pt} 2.5~  If 
                    $p_{+}$  satisfies   (\ref{lin2}), set $p_k=p_{+}$ and  go to Step 3.\\
\hspace*{44pt} 2.6~ Otherwise set $\gamma=\sigma\, \gamma$.\\
\hspace*{15pt} 3. Set $\gamma_k=\gamma$, $x_{k+1}=x_k+ p_k$.\\
\hspace*{15pt} 4. If $\|F_{k+1}\| =0$ stop.  \\   
\hspace*{15pt} 5. Choose $\beta_{k+1}$ such that $|\beta_{k+1}| \in  [\beta_{\min}, \beta_{\max}]$ .\\
}

The above results hold for any choice of the steplenght $\beta_k$ and Item 3 identifies  one occurrence
where the \name algorithm solves problem 
(\ref{ns}), i.e., $\{\|F_k\|\}$ converges to zero. In this section we complete the theoretical
 analysis of the \name algorithm by providing sufficient conditions that ensures
 that the sequence $\{\|F_k\|\}$ converges to zero. 
 
We start by recalling a simple result.
 \begin{lemma}
 Suppose that Assumption \ref{A2} holds. Then for $p_k = \pm  \gamma_k \beta_k F_k$, it holds
 \begin{equation}\label{eq:decF}
\small \|F_{k+1}\|^2 =  \left (1 \pm 2 \gamma_k \beta_k q((G_S)_k,F_k)  \pm 2 \frac{\gamma_k \beta_k}{\|F_k\|^2} \int_0^1 (F(x_k+p_k) - F(x_k))^T J(x_k+tp_k)F_k \, dt \right) \|F_k\|^2.
 \end{equation}
 \end{lemma}
 \begin{proof}
 Assume that $p_k = - \gamma_k \beta_k F_k$. Then, 
 \[
 \begin{array}{lcl}
 \|F_{k+1}\|^2&=&\|F_k\|^2+2\int_0^1 F(x_k+t p_k)^T J(x_k+tp_k) p_k\, dt  \\
 & =& \|F_k\|^2 - 2 \gamma_k \beta_k \int_0^1 F(x_k+tp_k)^T J(x_k+tp_k)F_k \, dt \\
 & =& \|F_k\|^2 - 2 \gamma_k \beta_k \int_0^1 F(x_k+tp_k)^T J(x_k+tp_k)F_k \, dt \\
&   &  \pm 2 \gamma_k \beta_k \int_0^1 F(x_k)^T J(x_k+tp_k)F_k \, dt \\
 & = & \|F_k\|^2 - 2 \gamma_k \beta_k F_k^T G_k F_k - 2 \gamma_k \beta_k \int_0^1 (F(x_k+p_k) - F(x_k))^T J(x_k+tp_k)F_k \, dt,\\
 \end{array}
 \]
 that gives (\ref{eq:decF}) using (\ref{symm}) and (\ref{RQ}). The case $p_k = + \gamma_k \beta_k F_k$ is analogous. 
 \end{proof}

Under specific assumptions on the Jacobian $J$, the following two theorems give conditions that ensures 
$F(x^*)=0$ where $x^*$ is the limit point  of $\{x_k\}$: Theorem  \ref{conv_glob_beta} concerns the cases when $J_S(x^*)$ is positive (negative) definite and when $J$ is symmetric too, Theorem  \ref{conv_glob_beta_ind}  regards the case when $J_S(x^*)$ is indefinite.

\begin{theorem}\label{conv_glob_beta}
Suppose that $F$ is continuously differentiable on $\IR^n$.
Let the positive sequence $\{\eta_k\}$ satisfy (\ref{etaserie}) and 
let $\{x_k\}$ be the sequence generated by the \name algorithm.
Moreover assume that  $J_S(x^*)$ is positive definite  at the limit point $x^*$ of $\{x_k\}$.
Letting $\sigmax(J(x^*))$ be the largest singular value of $J(x^*)$, if eventually  

    \begin{subequations}\label{eq:3}
    \noindent\begin{minipage}{0.4\textwidth}
\begin{equation}
\nu \ge \beta_k > \frac{\rho}{(1+\epsilon) \sigmax(J(x^*))} \label{cond_betarho}
\end{equation}
    \end{minipage}%
    \begin{minipage}{0.2\textwidth}\centering
    and
    \end{minipage}%
    \begin{minipage}{0.4\textwidth}
\begin{equation}
 \beta_k q((G_S)_k,F_k) > \frac 3 2 \rho,\label{cond_betaq}
\end{equation}
    \end{minipage}\vskip1em
\end{subequations}
\noindent
with $\rho\in (0,1)$ as in \req{lin1}-\req{lin2}  and for some $\epsilon\in (0,1)$ and $\nu>0$, then   $F(x^*)=0$. 
If  $\beta_k$  is either $\beta_{k,1}$ or $\beta_{k,2}$, only condition \req{cond_betaq} has to be satisfied to get $F(x^*)=0$. Moreover, for some  $\omega_1, \omega_2\in (0,1)$, sufficient conditions for  (\ref{cond_betaq})  to hold are
\begin{enumerate} 
\item 
if $\beta_k = \beta_{k,1}$ for $k$ large enough:
\begin{equation}\label{cond_gen22}
 \kappa(J_S(x^*))< \frac{2\omega_1}{3\rho};
\end{equation}
\item 
if $\beta_k = \beta_{k,2}$ for $k$ large enough:
\begin{equation}\label{cond_gen2}
 \kappa(J_S(x^*))< \omega_2 \sqrt{\frac{2}{3\rho}};
\end{equation}
\item if $J$ is symmetric and $\beta_k$ is either $\beta_{k,1}$ or $\beta_{k,2}$ for $k$ large enough:
\begin{equation}\label{cond_gen}
 \kappa(J(x^*))< \frac{2\omega_1}{3\rho};
\end{equation}
\end{enumerate}
where $\kappa(\cdot)$ is the 2-norm condition number.
\end{theorem}
\begin{proof}
Since $J_S(x^*)$ is assumed to be positive definite, continuity implies that there exists a scalar $\xi >0$  sufficiently small such that, for all 
$y \in \calB(x^*,\xi) = \{x \in \IR^n \ : \  \|x-x^*\|\le \xi\}$, 
$J_S(y)$ is positive definite and 
\begin{equation}\label{sigJS}
\lmin(J_S(y)) \geq (1-\epsilon)\lmin(J_S(x^*)), \mbox{ and } \lmax(J_S(y)) \leq (1+\epsilon)\lmax(J_S(x^*)),
\end{equation}
with $\epsilon \in (0,1)$. 
Moreover, the convergence of the sequence $\{x_k\}$ implies that
$x_{k-1}+tp_{k-1}$ and $x_k+tp_k$ both belong to $\calB(x^*,\xi)$ for large
enough $k$ and all $t \in [0,1]$.
As a consequence, reducing $\xi$ if necessary, we deduce that, for $k$ sufficiently large,
\begin{eqnarray*}
& \min\left[\lmin((G_S)_k),\lmin((G_S)_{k-1})\right] \geq (1-\epsilon)\lmin(J_S(x^*)), \label{sigGlbound}\\
& \max\left[\lmax((G_S)_k),\lmax((G_S)_{k-1})\right] \leq (1+\epsilon)\lmax(J_S(x^*)),\label{sigGubound}
\end{eqnarray*}
and by (\ref{eigenvalues}), 
\begin{equation}\label{eq:quoz2}
q((G_S)_k,F_k)  \in \left[
  \lmin((G_S)_k),\lmax((G_S)_k)\right]
\subseteq \left[
(1-\epsilon)\lmin(J_S(x^*)),
  (1+\epsilon)\lmax(J_S(x^*))\right].
\end{equation}
Finally, again by continuity,  reducing  $\xi >0$   if necessary, for all $y \in \calB(x^*,\xi)$ it holds 
\begin{equation}\label{sigJ}
\sigmax(J(y)) \leq (1+\epsilon)\sigmax(J(x^*)),\quad \sigmax(G_k) \leq (1+\epsilon)\sigmax(J(x^*)).
\end{equation}

Now, we consider (\ref{eq:decF}) and  $p_k=-\gamma_k \beta_k F_k$.
From the Mean Value Theorem  \cite[Lemma 4.1.9]{DS}, we have that
\[
\left |  \int_0^1 (F(x_k+tp_k)- F_k)^T J(x_k+tp_k)F_k \, dt \right | 
=
 \left |\int_0^1 \left(\int_0^1 J(x_k+\zeta\, tp_k)tp_k\,d\zeta\right)  J(x_k+tp_k)F_k \,dt \right|,
\]
$\zeta\in [0,1]$. Again, for $k$ sufficiently large, $x_k+\zeta \,tp_k\in \calB(x^*, \xi)$ for $t,\zeta\in [0,1]$.
Thus,  $p_k=-\gamma_k \beta_k F_k$ and (\ref{sigJ}) imply
\begin{eqnarray}
\left|  \int_0^1 (F(x_k+tp_k)- F_k)^T J(x_k+tp_k)F_k \, dt \right| 
&\le&\int_0^1 t\gamma_k   \beta_k \max_{z\in \calB(x^*, \xi)} \| J(z) \|^2   \|F_k\|^2  \, dt \nonumber  \\
&=& \frac{1}{2}\gamma_k    \beta_k  \max_{z\in \calB(x^*, \xi)} \sigmax(J(z))^2 \|F_k\|^2 \nonumber \\
&\le& \frac{1}{2}\gamma_k    \beta_k  (1+\epsilon)^2 \sigmax(J(x^*))^2 \|F_k\|^2. \nonumber
\end{eqnarray}
Combining this expression with (\ref{eq:decF}),  we have that for $k$ sufficiently large
\begin{eqnarray}
\|F_{k+1}\|^2& \le & \left (1 - 2 \gamma_k \beta_k q((G_S)_k,F_k)  + 2  \frac{\gamma_k \beta_k }{\|F_k\|^2}
\left|
 \int_0^1 (F(x_k+p_k) - F(x_k))^T J(x_k+tp_k)F_k \, dt \right | \right) \|F_k\|^2  \nonumber \\
             & \le & \left(1 - 2 \gamma_k \beta_k q((G_S)_k,F_k) +  \gamma_k^2 \beta_k^2 (1+\epsilon)^2 \sigmax(J(x^*))^2 \right) \|F_k\|^2.  \label{disbeta}
\end{eqnarray}

Thus, for $k$ sufficiently large, the linesearch condition \req{lin2} is satisfied  if  
$$
1 - 2 \gamma \beta_k q((G_S)_k,F_k)+ \gamma^2 \beta_k^2(1+\epsilon)^2 \sigmax(J(x^*))^2
\leq (1 - \rho \gamma)^2, 
$$
which is equivalent to
\begin{equation}\label{eq_lambda_bq}
\delta_2 \gamma^2 + 2 \delta_1 \gamma
\eqdef \left ( (1+\epsilon)^2\sigmax(J(x^*))^2 \beta_k^2 - \rho^2 \right )\gamma^2 + 2\left( \rho
-\beta_k q((G_S)_k,F_k)  \right)\gamma
\leq 0.
\end{equation}
Clearly \req{cond_betarho} implies that
 $ (1+\epsilon)^2\sigmax(J(x^*))^2 \nu^2 \ge \delta_2> 0$.
Moreover, if eventually  \req{cond_betaq} holds  then $\delta_1 <0 $ and (\ref{eq_lambda_bq}) is satisfied whenever
$\gamma \le \gamma^*= - 2\delta_1/\delta_2$. Now, $\gamma_*$ is uniformly bounded  below 
since $- \delta_1 \ge  \frac 1 2 \rho$, i.e., 
$\gamma^* \ge \frac {\rho}{\delta_2} \ge 
\bar \gamma \eqdef \rho/((1+\epsilon)^2\sigmax(J(x^*))^2 \nu^2)$.
Then, the mechanism of Step 3.6 of the \name algorithm guarantees that, for
$k$ sufficiently large, the loop in Step 2 terminates with $\gamma_k \geq
\min\{1, \sigma { \bar \gamma}\}$, and $ \bar \gamma $ independent of $k$. 
As a consequence, $\liminf_{k\rightarrow \infty}\gamma_k > 0$ and  by Item 2. in Theorem \ref{convF} we have that $F(x^*)=0$. 

We now show that when $\beta_k$ is either $\beta_{k,1}$ or $\beta_{k,2}$ for $k$ sufficiently large,
then only condition \req{cond_betaq} has to be satisfied to get $F(x^*)=0$. 

Let $\beta_k = \bbuno$. Using Item (ii) in  Lemma \ref{lem:bounds} and (\ref{bGspd}), we have that $\beta_k$ is positive and satisfies
\begin{equation}\label{b12gs}
\frac{1}{(1+\epsilon) \lmax(J_S(x^*))} \le \beta_{k} \le \frac{1}{(1-\epsilon) \lmin(J_S(x^*))}.
\end{equation}
By definition of $J_S$, $\|J_S(x^*)\|\le \|J(x^*)\|$, hence $\lmax(J_S(x^*)) \le \sigmax(J(x^*))$.
Therefore \req{cond_betarho} is satisfied  being $\rho \in (0,1)$ and setting $\nu =  1/((1-\epsilon) \lmin(J_S(x^*)))$.

Let $\beta_k=\bbdue$. Since  $\bbdue \le \bbuno$, the upper bound in \req{cond_betarho} is guaranteed from the discussion above. 
Moreover from \req{disbeta} and again from $\bbdue \le \bbuno$,
the linesearch condition \req{lin2} is satisfied  if  
\begin{equation}\label{eq_lambda_bq_b1}
\delta_2 \gamma^2 + 2 \delta_1 \gamma
\eqdef \left ( (1+\epsilon)^2\sigmax(J(x^*))^2 \beta_{1,k}^2 - \rho^2 \right )\gamma^2 + 2\left( \rho
-\beta_{2,k} q((G_S)_k,F_k)  \right)\gamma
\leq 0.
\end{equation}
Following the previous considerations on $\bbuno$, $\delta_2$ is positive. Further, 
using \req{cond_betaq} and  repeating the arguments above on 
the scalar $\gamma$ satisfying (\ref{eq_lambda_bq_b1}), the loop in Step 2  terminates with $\gamma_k \geq
\min\{1, \sigma { \bar \gamma}\}$, and $ \bar \gamma $ independent of $k$.

%
%

To conclude, as for Item 1., if $\beta_{k,1}$ is used eventually then \req{bGspd} and \req{eq:quoz2} give  $\beta_k q((G_S)_k,F_k)  
\ge  \frac{\omega_1}{ \kappa(J_S(x^*))}$
and trivially (\ref{cond_gen22}) implies (\ref{cond_betaq}) for all $k$ sufficiently large.

As for Item 2., if $\beta_{k,2}$ is used eventually then \req{bGspd2}, (\ref{sigJ}) and \req{eq:quoz2} give  $\beta_k q((G_S)_k,F_k)  
\ge  \frac{\omega_2^2}{ \kappa(J_S(x^*))^2}$ with 
$\omega_2= \frac{ (1-\epsilon  )\|J_S(x^*)\|}{(1+\epsilon  ) \|J(x^*)\|}$,
and  (\ref{cond_gen2}) implies (\ref{cond_betaq}) for all $k$ sufficiently large. 

Concerning Item 3.,  (\ref{cond_betaq}) reads $\beta_k q(G_k,F_k) > \frac 3 2 \rho$, and by 
Lemma \ref{lem:bounds} $\beta_{k,1}$ and $\beta_{k,2}$ are positive and  
$$
\beta_{k,1}\ge \beta_{k,2} \ge \frac{1}{\sigmax(G_{k-1})}\ge \frac{1}{(1+\epsilon)\sigmax(J(x^*))}. 
$$
Thus,  by  \req{eq:quoz2} it follows $\beta_k q(G_k,F_k)  \ge  \frac{\omega_1}{ \kappa(J(x^*))}$
and trivially (\ref{cond_gen}) implies (\ref{cond_betaq}) for all $k$ sufficiently large.
\end{proof}

\vskip 5pt
We remark that analogous conditions to \req{eq:3} can be derived for the case when 
$J_S(x^*)$ is negative definite. 
\begin{theorem}\label{conv_glob_beta_ind}
Suppose that $F$ is continuously differentiable on $\IR^n$.
Let the positive sequence $\{\eta_k\}$ satisfy (\ref{etaserie}) and 
let $\{x_k\}$ be the sequence generated by the \name algorithm. 
Moreover assume that  $J_S(x^*)$ is indefinite and $J(x^*)$ is nonsingular at the limit point $x^*$ of $\{x_k\}$.
If eventually
\begin{subequations}\label{eq:4}
    \noindent\begin{minipage}{0.42\textwidth}
\begin{equation}
\nu \ge |\beta_k | > \frac{\rho}{(1+\epsilon) \sigmax(J(x^*))}  \label{cond_betarho2}
\end{equation}
    \end{minipage}%
    \begin{minipage}{0.15\textwidth}\centering
    and
    \end{minipage}%
    \begin{minipage}{0.4\textwidth}
\begin{equation}
|\beta_k q((G_S)_k,F_k)| > \frac 3 2 \rho,\label{cond_betaq2}
\end{equation}
    \end{minipage}\vskip1em
\end{subequations}
\noindent
with $\rho\in (0,1)$ as in \req{lin1}-\req{lin2} and for some $\epsilon\in (0,1)$ and $\nu >0$,  then   $F(x^*)=0$. 
\end{theorem}
\begin{proof}
We observe that for $k$ sufficiently large, the inequalities \req{sigJS}-\req{eq:quoz2} hold  for some $\epsilon \in (0,1)$ .
Moreover, considering  $p_k = \pm \gamma_k \beta_k F_k$ and proceeding as in the proof of Theorem \ref{conv_glob_beta}, we get that for $k$ sufficiently large
the following inequality holds 
$$
\|F_{k+1}\|^2 \le \left(1 \pm 2 \gamma_k \beta_k q((G_S)_k,F_k) +  \gamma_k^2 \beta_k^2 (1+\epsilon)^2 \sigmax(J(x^*))^2 \right) \|F_k\|^2.
$$
Therefore the linesearch condition \req{lin2} is satisfied  if  
\begin{equation}\label{eq_lambda_bq2}
\delta_2 \gamma^2 + 2 \delta_1 \gamma
\eqdef \left ( (1+\epsilon)^2 \sigmax(J(x^*))^2 \beta_k^2 - \rho^2 \right )\gamma^2 + 2\left( \rho
\pm\beta_k q((G_S)_k,F_k)  \right)\gamma
\leq 0.
\end{equation}
Clearly  \req{cond_betarho2} implies that $ 
(1+\epsilon)^2 \sigmax(J(x^*))^2 \nu^2 \ge \delta_2>0$. 

We now show that \req{cond_betaq2} 
implies that $\delta_1>0$ so that we conclude that  $F(x^*)=0$ as in the proof of Theorem \ref{conv_glob_beta}.

Let us analyse the case $ \beta_k q((G_S)_k,F_k) < 0$ and consider the step $p_k=\gamma_k \beta_k F_k$.
Then condition \req{cond_betaq2}  means that $-\beta_k q((G_S)_k,F_k)  \ge \frac 3 2 \rho$, that is $\delta_1 = \rho + \beta_k q((G_S)_k,F_k)  < -\frac 1 2 \rho  <0$.
The case $\beta_k q((G_S)_k,F_k) > 0$ is analogous considering the step  $p_k= - \gamma_k \beta_k F_k$.
Now, repeating the arguments in Theorem \ref{conv_glob_beta} we conclude that  $\liminf_{k\rightarrow \infty}\gamma_k > 0$.
\end{proof}
} 
\section{Numerical experiments}
In view of our theoretical analysis and guidelines on steplength selection given in Section 4, we attempt
to tailor Barzilai and Borwein rules for unconstrained optimization
to spectral residual methods. In this section we  discuss 
several steplength rules for spectral residual methods and perform their experimental analysis using the \name algorithm described in 
Algorithm  \ref{pand_algo}.
Our test set consists of sequences of nonlinear systems arising in the solution of 
rail-wheel contact models and is described in details in Section \ref{treni}.


\name was implemented in Matlab (MATLAB R2019b) and  the experiments 
were carried out  on a Intel Core i7-9700K CPU @ 3.60GHz x 8, 16 GB RAM, 64-bit.

\subsection{Steplength rules} \label{sec::stepal}

We now present six  rules for the choice of the steplength in spectral residual methods  that were  used  
in our experiments. 
Besides the straightforward choice
of one of the two steplengths $\bbuno$, $\bbdue$, along all iterations, we consider adaptive strategies
that suitably combine them and parallel those 
used  for quadratic and nonlinear optimization problems. Below, given a scalar $\beta$,  
$T(\beta)$ is the thresholding rule which projects $|\beta|$ onto 
$\Ibeta \eqdef [ \beta_{\min}, \beta_{\max}]$ 
\begin{equation}\label{soglia}
    T(\beta) = \min\Big\{\beta_{\text{max}}, \max\big\{ \beta_{\text{min}},\big|\beta \big| \big\}\Big\}.
\end{equation}

\vskip 5pt
\begin{description}
\item{{\bf {\small{\rm {\bf   BB1}}} rule.}} By \cite{Sane,GS,Psane,Pand}, at each iteration let 
\begin{equation}\label{pbb1}
\beta_k=
\begin{cases}
\bbuno & \text{if  }  \ |\bbuno|\in \Ibeta \\
T(\bbuno)  & \text{otherwise} 
\end{cases}    
\end{equation}
\vskip 3pt
\item{{\bf {\small{\rm {\bf   BB2}}} rule.}} At each iteration let 
\begin{equation}\label{pbb2}
\beta_k=
\begin{cases}
\bbdue & \text{if  }  \ |\bbdue|\in \Ibeta  \\
T(\bbdue)  & \text{otherwise} 
\end{cases}   
\end{equation}
\vskip 3pt
\item{{\bf {\small{\rm {\bf   ALT}}} rule.}} Following 
\cite{DF, GS}, at each iteration let us alternate  between $\bbuno$ and $\bbdue$:
\begin{eqnarray}
& & 
\beta^{{\small{\rm ALT}}}_k= 
\begin{cases}
\bbuno & \text{for } k \mbox{ odd} \\
\bbdue & \text{otherwise} 
\end{cases}   \label{palt1} \\
& & \nonumber \\
& & \beta_k= 
\begin{cases}
\beta^{{\small{\rm ALT}}}_k & \quad \text{if} \ \ \ |\beta_k^{{\small{\rm ALT}}}| \in \Ibeta \\
\bbuno                      & \quad \text{if $k$  even,} \ \ |\bbuno| \in \Ibeta,\ |\bbdue| \notin \Ibeta  \\
\bbdue                      & \quad \text{if $k$  odd,} \ \  \ |\bbdue| \in \Ibeta ,\ |\bbuno| \notin \Ibeta\label{palt2}  \\
T(\beta^{{\small{\rm ALT}}}_k) & \quad \text{otherwise} 
\end{cases}
\end{eqnarray}

\vskip 3pt
\item{{\bf {\small{\rm {\bf   ABB}}} rule.}} Following \cite{Zhou}  and {\small{{\rm ABB}}} rule in \cite{FZanni}, 
we define the  Adaptive Barzilai-Borwein ({\small{\rm ABB}}) rule  as follows. Given $\tau\in (0,1)$, let
\begin{eqnarray}
& &  \beta^{{\small{\rm ABB}}}_k(\xi_1,\xi_2) =
\begin{cases}       
\xi_2 & \text{if} \quad \displaystyle  \frac{\xi_2}{\xi_1} < \tau \\
\xi_1 & \text{otherwise} 
\end{cases}\label{pabb1}
\end{eqnarray}
for some given $\xi_1, \, \xi_2$.  Then
\begin{eqnarray}
& & \beta_k= 
\begin{cases} 
\beta^{{\small{\rm ABB}}}_k(\bbuno, \bbdue)  & \quad \text{if } \  \ |\bbuno|,|\bbdue| \in \Ibeta\\
\bbuno                      & \quad \text{if } \ \ |\bbuno| \in \Ibeta ,\ |\bbdue| \notin \Ibeta\\
\bbdue                      & \quad \text{if } \ \ |\bbdue| \in \Ibeta ,\ |\bbuno| \notin \Ibeta \\
\beta^{{\small{\rm ABB}}}_k(T(\bbuno), T(\bbdue))  & \quad \text{otherwise} \label{pabb2}
 \end{cases}
\end{eqnarray}
Observe that  a large value of $\tau$ promotes the use of $\bbdue$ with respect to $\bbuno$. The rule 
allows to switch between the steplengths $\bbuno$ and $\bbdue$ and was originally motivated by the behaviour of the Barziali and Borwein method applied to convex and quadratic minimization 
problem (see \cite{Zhou, FZanni} and our discussion below Lemma \ref{prop::sym}).

\vskip 3pt
\item{{\bf {\small{\rm {\bf ABBm}}} rule.}} This rule elaborates the {\small{{\rm ABBminmin}}} rule given in \cite{FZanni}, 
taking into account that $\bbdue$ may be  negative along iterations. 
Let $m$ be a  nonnegative integer,  and 
\begin{equation}\label{betamem}
\begin{array}{l}
\widetilde{\beta}_{k,2} = 
\begin{cases}
\bbdue  & \text{if} \ \ \ |\bbdue| \in  \Ibeta \\
T(\bbdue) & \text{otherwise} 
\end{cases}\\
\\
j^* = \argmin\{|\widetilde{\beta}_{j,2}|: j=\max\{1, k-m\}, \dots, k\}.
\end{array}
\end{equation}
Given $\tau \in (0,1)$, we fix $\beta_k$ as follows
\begin{eqnarray} 
& & \beta^{{\small{\rm ABBm }}}_k(\xi_1,\xi_2)= 
\begin{cases}
\widetilde{\beta}_{j^*,2} & \text{if} \quad  \displaystyle\frac{\xi_2}{\xi_1} < \tau \\
\xi_1& \text{otherwise}
\end{cases} \label{pabbmin1}\\
& & \nonumber \\
& & \beta_k= 
\begin{cases} 
\beta^{{\small{\rm ABBm}}}_k(\bbuno,\bbdue)  & \quad \text{if } \ \ |\bbuno|,|\bbdue| \in \Ibeta \\
\bbuno                      & \quad \text{if } \ \ |\bbuno| \in \Ibeta  ,\ |\bbdue| \notin \Ibeta \\
\bbdue                      & \quad \text{if } \ \ |\bbdue| \in \Ibeta  ,\ |\bbuno| \notin \Ibeta \\
\beta^{{\small{\rm ABBm}}}_k(T(\bbuno), T(\bbdue)) 
 & \quad \text{otherwise}    \label{pabbmin2}
 \end{cases}
\end{eqnarray}

Again, a large value of $\tau$ promotes the use of a step from BB2 rule instead of $\bbuno$. 
In case $|\bbuno|,|\bbdue| \in \Ibeta$ and $\displaystyle\frac{\bbdue}{\bbuno} < \tau$, the smallest absolute value 
 $\widetilde{\beta}_{j^*,2}$  over the last $m+1$ iterations is selected; {\color{black}
taking into account that $\widetilde{\beta}_{j,2}$ for $j=\max\{1, k-m\}, \dots, k$ can be negative,  the rationale for selecting $\widetilde \beta_{j^*,2}$ in (\ref{pabbmin1}) is to mitigate the nonmonotone behavior of the objective function \cite{FZanni}. 
Consequently,  smaller steplengths are expected using the {\small{\rm  ABBm}} rule than using  the {\small{\rm  ABB}} rule.}
\vskip 3pt
\item{{\bf {\small{\rm {\bf   DABBm}}} rule.}}
Following \cite{Bon, Zanni3},  a dynamic threshold $\tau_k\in (0,1)$ can be used in place of the prefixed threshold $\tau$ in \req{pabbmin1}. 
Given $\widetilde \beta_{k,2}$ and $j^*$ in \req{betamem}, we  propose the  rule defined as 
\begin{eqnarray} 
& & \beta^{{\small{\rm DABBm }}}_k(\xi_1,\xi_2)= 
\begin{cases}
\widetilde{\beta}_{j^*,2} & \text{if} \quad  \displaystyle\frac{\xi_2}{\xi_1} < \tau_k \\
\xi_1& \text{otherwise}
\end{cases} \label{pdabbmin1}\\
& & \nonumber \\
& & \beta_k= 
\begin{cases} 
\beta^{{\small{\rm DABBm}}}_k(\bbuno,\bbdue)  & \quad \text{if } \ \ |\bbuno|,|\bbdue| \in \Ibeta\\
\bbuno                      & \quad \text{if } \ \ |\bbuno| \in \Ibeta,\ |\bbdue| \notin \Ibeta\\
\bbdue                      & \quad \text{if } \ \ |\bbdue| \in \Ibeta ,\ |\bbuno| \notin \Ibeta \\
\beta^{{\small{\rm DABBm}}}_k(T(\bbuno), T(\bbdue)) 
 & \quad \text{otherwise}    \label{pdabbmin2}
 \end{cases}
\end{eqnarray}
with the dynamic threshold set as 
\begin{eqnarray}
& & \tau_k = \min \Big\{ \tau, \|F_k\|^{ 1/(2+{b_t}^2) } \Big\}, \label{tauk2}\\
& & 
b_t = \max \{b_j : j=\max\{1, k-w\}, \dots, k\}\label{tauk}.
\end{eqnarray}
Here $\tau \in (0,1)$ is an upper bound on the value of $\tau_k$, $w$ is a nonnegative integer and
$b_j$ denotes the number of backtracks performed at iteration $j$ (see Step 2 of Algorithm \ref{pand_algo}). 
If $\|F_k\|$ is getting small and  
the number of performed backtracks in the last $w+1$ iterations is small, then  (\ref{tauk2}) 
promotes the use of steplength from  BB1 rule, i.e.,  larger steplengths which can speed convergence to a zero  of $F$.
On the other hand, when the number of backtracks performed along previous iterations is large
and $\tau$ is large, 
the use of the smaller steplength from BB2 rule is encouraged.
\end{description}
{\color{black} We conclude the discussion on steplenght selection, noting that 
conditions \req{eq:3} and \req{eq:4} for the convergence of $\{x_k\}$ to a zero of $F$ apply to all
our rules.}
\vskip 5pt
The rules and parameters used in our experiments are summarized in Table \ref{tab:pandv}.
\begin{table}[h]
\small
\centering
\begin{tabular}{ l|l }
\toprule 
\multicolumn{ 1}{c|}{Rule}   & \multicolumn{ 1}{c}{   $ \beta_k$} \\
\midrule
\BBu      & $\beta_k$ in (\ref{pbb1}) \\
\BBd    & $ \beta_k$ in (\ref{pbb2}) \\
\BBalt     & $ \beta_k$ in (\ref{palt1}), (\ref{palt2}) \\
\abbu    & $ \beta_k$ in (\ref{pabb1}), (\ref{pabb2}) with $\tau=0.1$ \\
\abbo   & $\beta_k$ in  (\ref{pabb1}), (\ref{pabb2}) with $\tau=0.8$ \\
\abbminu & $ \beta_k$ in (\ref{betamem})-(\ref{pabbmin2}) with $\tau=0.1$, $m=5$\\
\abbmino &  $ \beta_k$ in (\ref{betamem})-(\ref{pabbmin2}) with $\tau=0.8$,  $m=5$\\
\dabbmino  &  $ \beta_k$ in (\ref{betamem}), (\ref{pdabbmin1})-(\ref{tauk}) with $\tau=0.8$, $m=5$, $w=20$\\
\bottomrule
\end{tabular}  
\caption{Steplength's rules in   \name implementation.}\label{tab:pandv}
\end{table}

\subsection{Problem set: nonlinear systems arising from rolling contact models}\label{treni}
Rolling contact is a fundamental issue in mechanical engineering and plays a central role in many important applications such as
rolling bearings and wheel-rail interaction \cite{Kalker1, Kalker2}. In order to perform simulations  of 
complex mechanical systems with a good tradeoff between accuracy and efficiency, three working hypotheses are usually made
in modelling rolling contact: non-conformal contact, i.e., the typical dimensions of the contact area are negligible 
if compared to the curvature radii of the contact body surfaces; planar contact, i.e., the contact area is contained in a plane; 
half-space contact, i.e., locally, the contact bodies are viewed as three-dimensional half-spaces \cite{Kalker1, Kalker2}. 
In this framework, we focus on the Kalker's rolling contact model which represents a relevant and general model in contact mechanics.

The solution of Kalker's rolling contact model can be performed using different approaches.  
The approach in \cite{Voll, Voll2} calls for the solution of constrained optimization problems while  
the so-called {\small{\rm CONTACT}} algorithm \cite{Kalker2} gives rise to sequences of nonlinear systems.
Our problem set derives from the application of {\small{\rm CONTACT}} algorithm; here we describe in which phase of the Kalker's model solution they arise and give some of their features. We refer  to Appendix \ref{sec:model} for a sketch of Kalker's model, its discretization, and the Kalker's {\small{\rm CONTACT}} algorithm. 

Kalker's {\small{\rm CONTACT}} algorithm determines the normal pressure, the tangential pressure, the contact area, the adhesion area and the sliding area
in the contact between two elastic bodies
and relies on the elastic decoupling between the normal contact problem and the tangential contact problem. 
Such problems are solved separately; first
the normal problem is solved via the the so-called {\small{\rm NORM}} algorithm, second the tangential 
problem is solved via the so-called {\small{\rm TANG}} algorithm. 
Algorithms {\small{\rm NORM}} and {\small{\rm TANG}} are expected to identify the  
elements in the contact area and in the adhesion-sliding areas, respectively. 
These algorithms are applied sequentially and repeatedly until the values of the computed pressures undergo a sufficiently small change that 
suggests their reliable approximation; in general, a few repetitions of {\small{\rm NORM}} and {\small{\rm TANG}} algorithms are required.
Each repetition of {\small{\rm NORM}} algorithm calls for the solution of a sequence of linear systems while each 
repetition of {\small{\rm TANG}} algorithm calls for the solution of a sequence of linear and nonlinear systems.
Computationally, the major bottleneck is the numerical solution of the sequence 
of nonlinear systems generated in the {\small{\rm TANG}} phase.
Importantly, each {\small{\rm CONTACT}} iteration requires few repetitions of {\small{\rm TANG}} algorithm but the {\small{\rm CONTACT}} algorithm is performed for several time instances\footnote{In Appendix \ref{sec:model} see: (\ref{contact_problem}) for the form of  normal contact problem and 
tangential contact problem, (\ref{nleq_d}) for the form of the nonlinear systems to be solved, Figure \ref{fig:archkca} for the flow of  Kalker's {\small{\rm CONTACT}} algorithm.}.

\begin{figure}[htbp]
	\centering
  \includegraphics[width=.8\textwidth, height=.6\textheight, keepaspectratio]{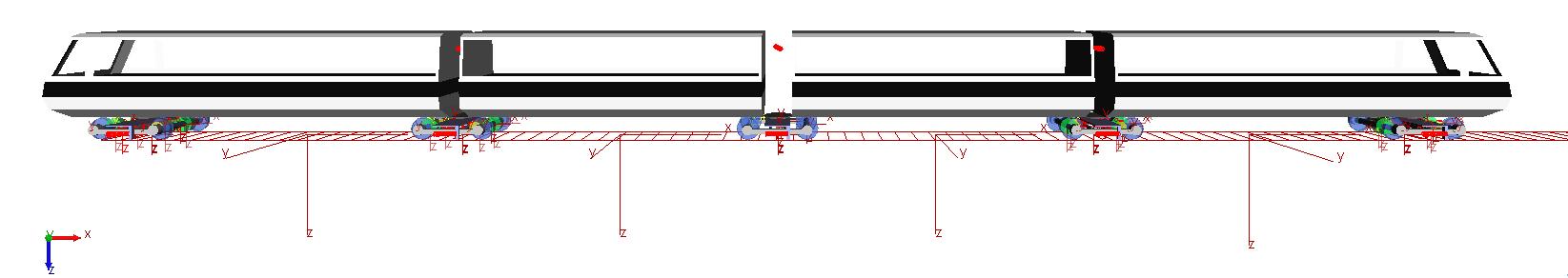}
	\caption{Multibody model of the benchmark vehicle.}
	\label{fig:treno}
\end{figure}

Our tests were made on wheel-rail contact in railway systems. 
The benchmark vehicle is a driverless subway vehicle, designed by Hitachi Rail 
on MLA platform (Light Automatic Metro).
The vehicle is a fixed-length train composed of four carbodies and five bogies
(four motorized and one, the third, trailer), see Figure \ref{fig:treno}.
The multibody model has been realized in the Simpack Rail environment \cite{Simpack}.
We considered a train route of length $400m$ including a typical
railway curved track
characterized by three significant parts: two straight lines (from $0 m$ to $70 m$ and from $233 m$ to $400 m$), the curve (from $116 m$ to $186 m$) and two cycloids (from $70 m$ to $116 m$ and from $186 m$ to $233 m$) which smoothly connect the straight lines and the curve
in terms of curvature radius. The radius of the curve is $500 m$. In this analysis, we focused on the contact between the first vehicle wheel and the rail; since the vehicle length is equal to $45.7 m$, at the beginning of the dynamic simulation the considered wheel starts in the position $45.7 m$ along the track.
We performed a simulation in an interval of  10 seconds using 500 time steps, 
which amounts to 500 calls to {\small{\rm CONTACT}} algorithm, 
for  train speeds with magnitude $v$ taking the values: $ {v}=10\: m/s$ and $ {v}=16\: m/s$. Accordingly, during the whole simulation the considered wheel travels along the track a distance equal to $100 m$ and $160 m$, respectively. The traveling velocities considered give a realistic lateral acceleration along the curve according to the current regulation in force in the railway field.  

Two sets of experiments were performed\footnote{The data that support the findings of this study are available from the corresponding author upon reasonable request.}. 
First, we solved a large number of sequences of nonlinear systems arising from 
wheel-rail contact in railway systems by 
the eight \name variants based on   the rules in  Table \ref{tab:pandv}. Second, we compared experimentally
the best performing \name variant and  a standard Newton trust-region  when embedded in the {\small{\rm CONTACT}} algorithm.

The set of test problems used in the first part of the experiments  was generated implementing 
the {\small{\rm CONTACT}} algorithm in Matlab  and using
a  standard trust-region Newton method\footnote{The code in \cite{Tresnei} was applied using the default setting and dropping bound constraints on the unknown.} for solving 
the arising nonlinear systems. Afterwards, a representative subset of the nonlinear systems was selected 
to form our problem set. Specifically,  six sequences of nonlinear systems generated   
by the {\small{\rm CONTACT}} algorithm and corresponding to six consecutive time instances for
each track section (straight line, cycloid and curve) and for each velocity were selected. 
Such sequences are representative of the systems arising throughout the whole simulation and 
allow a fair analysis of \name on nonlinear systems from  a real application. 
 Table \ref{stretch} summarizes the features of the sequences: magnitude of the
train velocity  $v$, section of the route, time instances, number of nonlinear systems
in the sequence, dimension $n$ of the systems (proportional to the number of mesh nodes in the potential contact area). 
A typical feature of the contact model is that $n$ increases as the velocity increases and when the train curves along the route (i.e., the track curvature increases). The total number of systems associated to ${v}=10\: m/s$ and ${v}=16\: m/s$ is 121 and 153 respectively.
\begin{table}[h!]
    \centering
\begin{tabular}{clcc c}
\toprule
${v}$($m/s$)   &   Track Section     &  Time Instances   & Number of Systems & $n$  \\
\midrule
    &   Straight line    &  100-105   &  10             &  156             \\
10   &   Cycloid    &   300-305   &    56     &     897           \\
     &   Curve     &     450-455   &     55     &       1394       \\
\midrule

    &    Straight line   &   50-55    &  8   &   156    \\
16    &  Cycloid     &   150-155  & 63   &   1120       \\
     &  Curve    &   350-355    &   82   &    1394           \\
\bottomrule
\end{tabular}
\caption{Sequences of nonlinear systems forming the first problem set.}\label{stretch}
\end{table}

\subsection{Numerical results}\label{numer}
 In this section we present the performance of \name algorithm. 
 The results presented concern the solution of the sequences of nonlinear systems summarized in Table \ref{stretch} and a comparison between
the best performing \name variant and  a standard Newton trust-region  method when embedded in the {\small{\rm CONTACT}} algorithm.

\name algorithm was   implemented as described in Section \ref{sec::stepal} and 
with  parameters 
$$\beta_{\text{min}}= 10^{-10},\ \  \beta_{\text{max}}=10^{10},\ \  \rho= 10^{-4}, \ \ \sigma=0.5,\ \ \eta_k = 0.99^k(100+\|F_0\|^2)  \  \forall k\ge 0,$$
see \cite{Pand}. The null vector $x_0=0$  was chosen as initial guess.
A maximum number of iterations   and $F$-evaluations  equal to $10^5$ was imposed and a maximum number of backtracks 
equal to 40 was allowed at each iteration.  
The procedure was declared successful when 
\begin{equation}\label{stop}
\| F_k  \| \le 10^{-6}.    
\end{equation}
A failure was declared either because the assigned maximum number of iterations or $F$-evaluations or backtracks is reached,
or because $\|F\|$ was not reduced for 50  consecutive  iterations.

We now compare the performance of all  the variants of \name method in the solution of the sequences of nonlinear systems 
in Table \ref{stretch}. Further, in light of the theoretical investigation presented in this work, we analyze in details  the results obtained with  \BBu and \BBd rule and support the use of rules that switch between the two steplengths.

 \begin{figure}[h!]
\centering
   {\includegraphics[width=.98\textwidth]{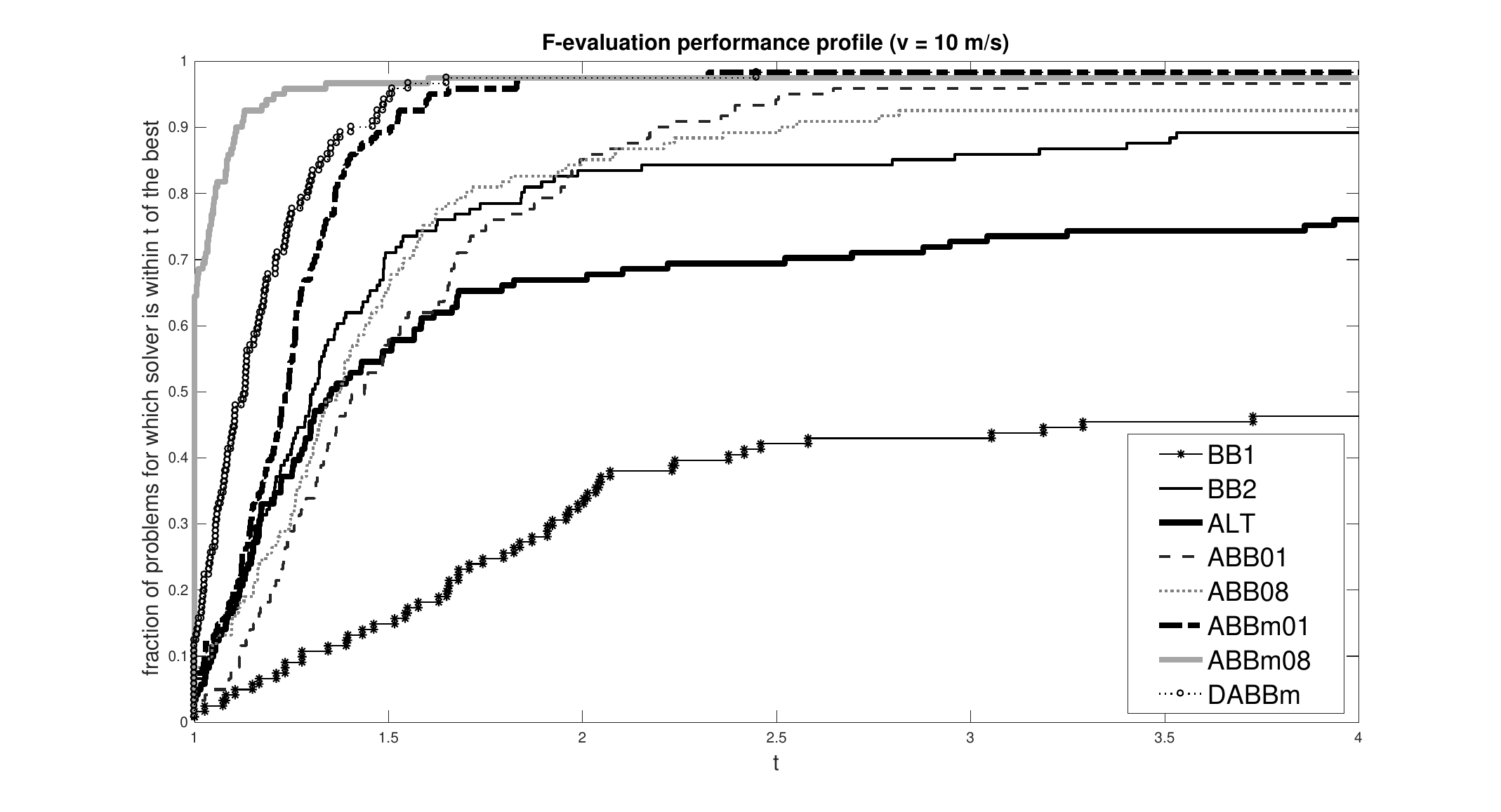}} \quad
  {\includegraphics[width=.98\textwidth]{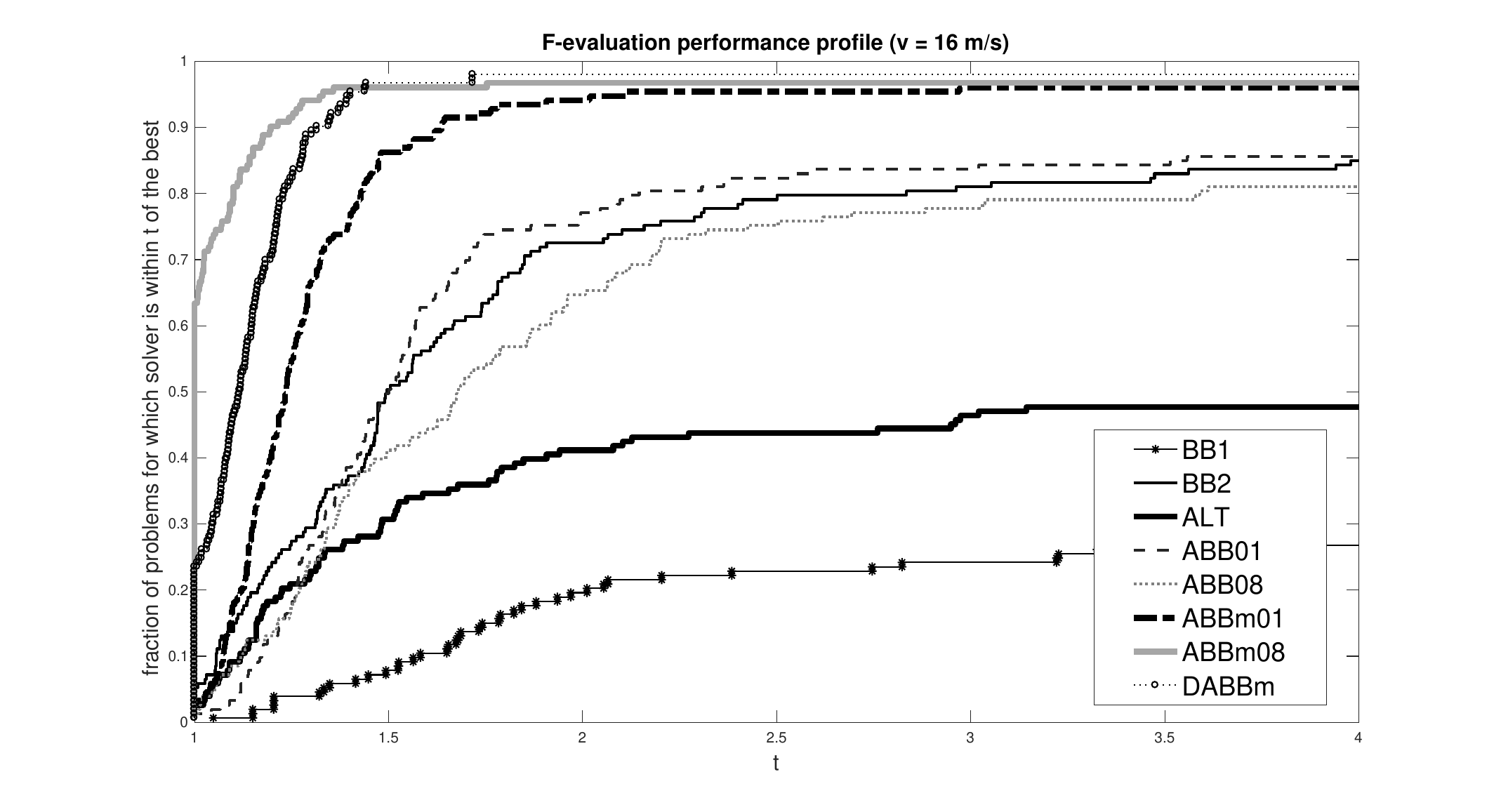}} 
 \caption{$F$-evaluation performance profiles of \name method. Upper: $ {v}=10\: m/s$, Lower: $ {v}=16\: m/s$.}
 \label{Perf}
 \end{figure}
 
Figure \ref{Perf} shows the performance profiles  \cite{dm} 
in terms of $F$-evaluations employed by the \name variants  
for solving the sequence of systems generated both with  ${v}=10\: m/s$ (121 systems) (upper) and with 
${v}=16\: m/s$ (153 systems) (lower) and highlights that the choice of the steplength is crucial for both efficiency and robustness. The complete results are reported in Appendix \ref{sec::fullres}.
We start observing that   \BBd rule outperformed {\small{\rm BB1}} rule; 
in fact  the latter shows the  worst behaviour both in terms of 
efficiency and in terms  of number of systems solved.
Alternating $\bbuno$ and $\bbdue$ in \BBalt rule without taking into account the magnitude 
of the two scalars improves performance over \BBu rule but is not competitive with {\small{\rm BB2}} rule.
On the other hand, 
the variants of \name using adaptive strategies are the most robust, i.e., they solve the largest number of problems,
and efficient. 
Specifically,  comparing {\small{\rm ABB}}, {\small{\rm ABBm}} and  {\small{\rm DABBm}} rules, the most effective steplength selections are {\small{\rm ABBm}} and  {\small{\rm DABBm}}. Using \abbminu rule, 98.3\% (2 failures)  and 96.1\% (6 failures) out of the total  number of systems were solved successfully
for ${v}=10\: m/s$ and ${v}=16\: m/s$
respectively; using \abbmino rule,   98.3\% (2 failures)  and 96.7\% (5 failures) of the total  number of systems
were solved successfully with ${v}=10\: m/s$ and ${v}=16\: m/s$
respectively; using the dynamic selection {\small{\rm DABBm}},  the largest number of systems was solved successfully, i.e., 
99.2\% (1 failure)  and 98\% (3 failures) out the total  number of systems with ${v}=10\: m/s$ and ${v}=16\: m/s$
respectively.
Overall, \abbmino rule gives rise to the most efficient algorithm for both velocity values and the profile related to \BBd rule is within a factor 2 of it  in roughly the 80\% and the 70\% of the runs for  ${v}=10\: m/s$ and ${v}=16\: m/s$, respectively.

Let us now focus on the performance  \name coupled with  \BBu and \BBd rules. 
As a representative run of our numerical experience reported in Appendix \ref{sec::fullres},  we consider 
the  nonlinear system arising with ${v}=16\: m/s$,
at time $t=150$,  iteration 2 of the {\small{\rm CONTACT}} algorithm and iteration 2 of the 
{\small{\rm TANG}}  algorithm (system 150\_2\_2 in Table \ref{tab:v16_cy}).
\begin{figure}[h!]
\centering
   {\includegraphics[width=\textwidth]{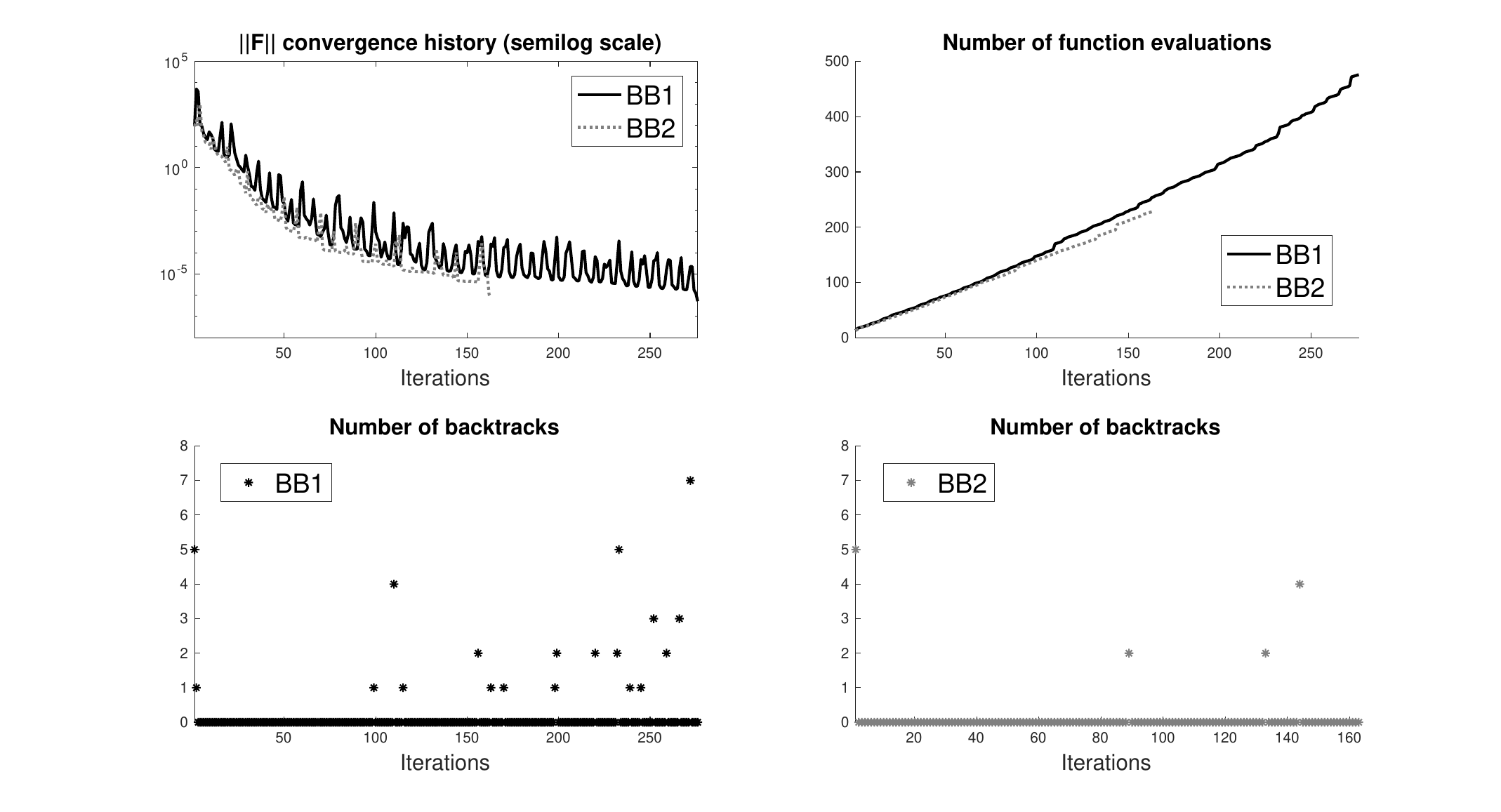}} \quad
 \caption{\name with \BBu rule vs \name with \BBd rule  on a single nonlinear system.}
  \label{BB1BB2}
 \end{figure}
In the upper part  of Figure \ref{BB1BB2} we display  $\|F\|$ along iterations and the number of $F$-evaluations performed. 
We note that using the stepsize $\bbuno$ causes a highly nonmonotone  behavior of $\|F\|$ and  such behaviour is not productive for
convergence;  using \BBu rule 276 iterations and 476 $F$-evaluations are performed while using \BBd  rule 163 iterations and 228 $F$-evaluations are required.
The distinguishing feature of these runs is the high  number of backtracks performed using $\bbuno$
at some iterations,
as reported at the bottom part of the figure where the number of backtracks versus iterations is reported for both \name variants.
This behaviour is in accordance with the analysis in Section \ref{sec::stepnm}: since $\bbuno$ can be arbitrarily larger 
than $\bbdue$ in the indefinite case, the need to perform a large number of backtracks to enforce approximate norm decrease is likely to occur 
in case $\bbuno$ is taken as the initial steplength.
Such observation supports the  use of $\bbdue$; the benefit from using 
shorter steps is  further shown by the performance of {\small{\rm ABBm}} over {\small{\rm ABB}}, 
the former  tends to take shorter steps than the latter by exploiting  the iteration history and  
results to be more effective. 

We conclude our experimental analysis using a spectral residual method 
in the {\small{\rm CONTACT}} algorithm. 
To this purpose, we compare two implementations of {\small{\rm CONTACT}} algorithm which differ only 
in the nonlinear solver for the nonlinear systems arising in the {\small{\rm TANG}} algorithm.
The first implementation ({\small{\rm CONTACT-NTR}}) uses a standard Newton trust-region method  and the
second one ({\small{\rm CONTACT-DABBm}})
uses {\small{\rm DABBm}} which turned out to be the more robust  \name version  
in the analysis above (see Figure \ref{Perf}).  
As a standard Newton trust-region method, we used the Matlab code  proposed in \cite{Tresnei}; default parameters were used and bound constraints on the unknown were dropped using the setting indicated in the code. The Jacobian matrix of $F$ was approximated by finite differences.

\begin{figure}[h!]
\centering
  {\includegraphics[width=\textwidth]{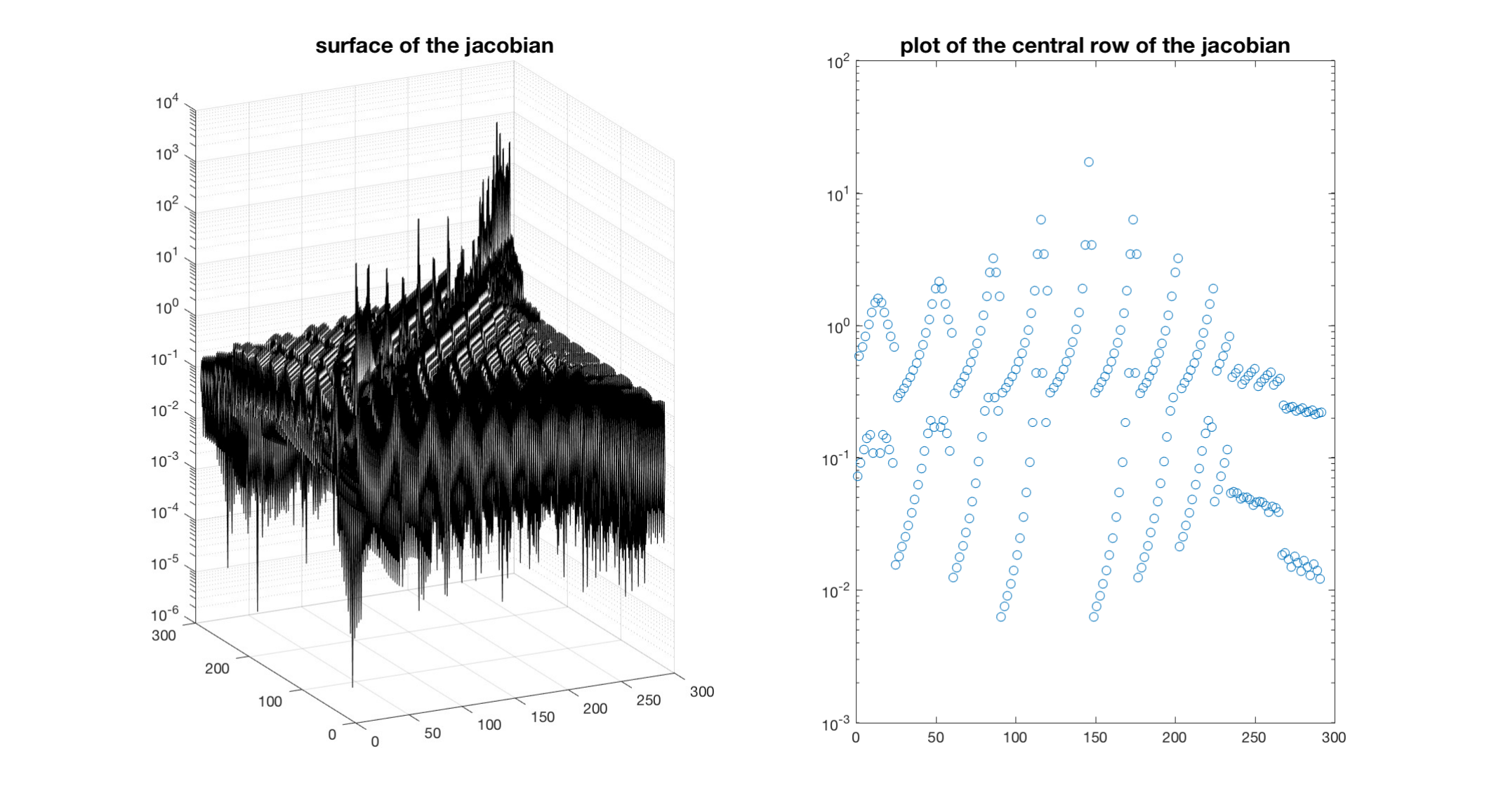}} 
 \caption{Jacobian matrix: surface of the full matrix and plot of the central row (base 10 logarithm of the absolute values).}
 \label{surf}
 \end{figure}

As a preliminary issue, we observe that the Jacobian matrices of $F$ are dense through the iterations;  thus they cannot be formed as a low computational cost by finite difference procedures for sparse matrices \cite{cpr}. 
We also observed in the experiments that the Jacobian matrices are nonsymmetric, do not have dominant diagonals and 
they are not close to diagonal matrices.  For example, let us  consider the Jacobian matrix of the system corresponding to speed $ {v}=16\: m/s$,
curve track section, instant $t=355$, iteration 2 of the {\small{\rm CONTACT}} and iteration 4 of the {\small{\rm TANG}} algorithm
(355\_2\_4 in Table \ref{tab:v16_cur}). It has dimension $292\times 292$ and, evaluated at the final iterate computed using \abbmino rule,
$96.18 \%$ of  its elements are nonzero. The structure of the Jacobian can be observed in Figure \ref{surf} where the absolute values of its elements are plotted in a logarithmic scale
(the surface of the full matrix on the left and a plot of the row 146 on the right). This structure is observed along all the iterations of the nonlinear system solvers and is common to all sequences generated by the {\small{\rm CONTACT}} algorithm.

In our implementation, {\small{\rm CONTACT}} algorithm terminated when the relative error between two successive values of the computed pressures dropped below  $10^{-4}$ or a maximum of 20 alternating cycles between {\small{\rm NORM}} and {\small{\rm TANG}} was reached.
Both nonlinear solvers were run until the stopping rule (\ref{stop}) is met.
We ran {\small{\rm CONTACT-NTR}} and {\small{\rm CONTACT-DABBm}}
over the whole track for both velocities, that is we considered the whole sequence of 500 time steps. {\small{\rm CONTACT-NTR}}  generated  3759 and 5353 nonlinear systems for $ {v}=10\: m/s$ and $ {v}=16\: m/s$, respectively and {\small{\rm CONTACT-DABBm}} generated
4496 and 5494 nonlinear systems for the two velocities.

\begin{figure}[h!]
\centering
   {\includegraphics[width=\textwidth]{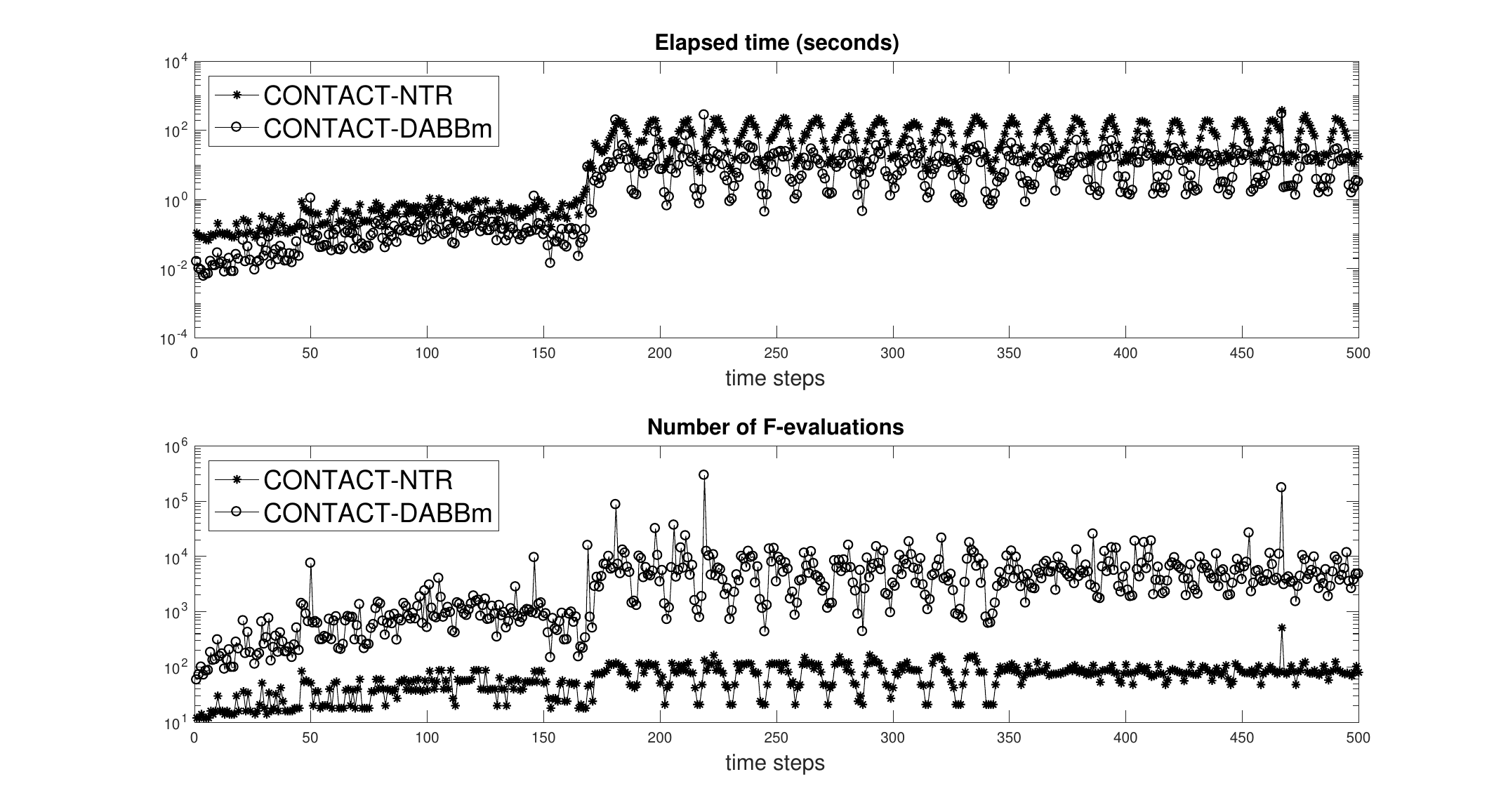}} \quad
 \caption{Comparison between  {\small{\rm CONTACT-DABBm}} and 
  {\small{\rm CONTACT-NTR}}, $ {v}=10\: m/s$: number of $F$evaluations and elapsed time in seconds  (logarithmic scale).}
 \label{FVtime}
 \end{figure}

\begin{figure}[h!]
\centering
  {\includegraphics[width=\textwidth]{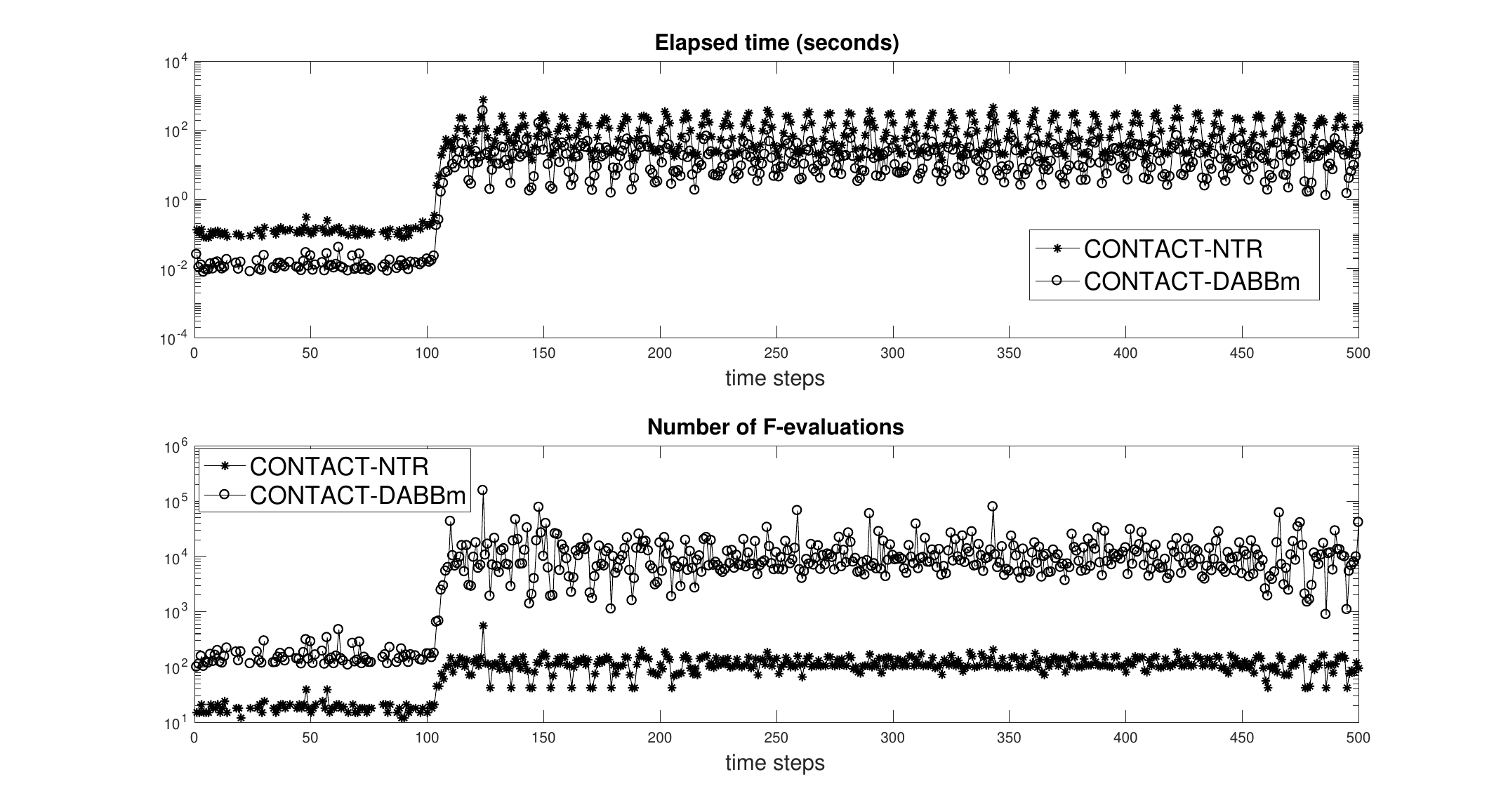}} \quad
 \caption{Comparison between   {\small{\rm CONTACT-DABBm}} and 
  {\small{\rm CONTACT-NTR}},  $ {v}=16\: m/s$: number of $F$evaluations and elapsed time in seconds (logarithmic scale).}
 \label{FVtime2}
 \end{figure}

As a first remark, both procedures successfully solved the contact model described above and were reliable and accurate 
in the  numerical simulation of wheel-rail interaction.
Secondly, the use of the spectral residual method
yields a gain  in terms of time with respect to 
the use of a standard Newton method where finite difference approximation
of Jacobian matrices is employed; this feature derives from the fact that
spectral residual method is derivative-free and does not ask for the solution of linear systems.
Figures \ref{FVtime} and \ref{FVtime2} show the comparison of the two {\small{\rm CONTACT}} implementations  in terms of number of $F$-evaluations (excluding those needed to approximate the Jacobian matrices) and execution elapsed time. From the plots we observe  that {\small{\rm CONTACT-DABBm}} takes a  larger number of $F$-evaluations than {\small{\rm CONTACT-NTR}} but it is faster. 
Over the whole time interval, {\small{\rm CONTACT-DABBm}} employs 1 hour, 19 mins and 2 hours, 28 mins to solve the generated nonlinear systems with $ {v}=10\: m/s$ and $ {v}=16\: m/s$, while {\small{\rm CONTACT-NTR}} takes 7 hours and 49 mins 
and 12 hours and 41 mins, respectively.

 \section{Conclusions}
 The numerical behaviour of spectral residual methods for nonlinear systems strictly depends  on the choice of the spectral steplength. Although most of the works on this subject make use of the
 stepsize $\bbuno$, known results on the spectral gradient methods for unconstrained optimization suggest
that a suitable combination of the stepsizes $\bbuno$ and $\bbdue$ could be of benefit for spectral residual methods as well. 
This work aims to contribute to this study by providing a first systematic analysis of the  stepsizes $\bbuno$ and $\bbdue$.  Moreover, practical guidelines for  dynamic choices of the steplength are derived from new theoretical results in order to increase both the robustness and the efficiency of spectral residual methods. Such findings have been extensively tested and validated on sequences of nonlinear systems arising in the solution of a contact wheel-rail model.

\subsection*{\footnotesize Acknowledgments}
INdAM-GNCS partially supported the second, the third and the fourth author under Progetti di Ricerca 2019 and 2020.\\

\noindent {\bf Declarations}\\

\noindent 
{\bf Conflict of interest} The authors declare that they have no conflict of interest.\\

\noindent 
{\bf Funding} Open access funding provided by Universit\`a di Bologna within the CRUI-CARE
Agreement.
\appendix


\section{ Kalker's contact model and  {\small{\bf {\rm CONTACT}}}  algorithm}
\label{sec:model}

\begin{figure}[htbp]
	\centering
	  \includegraphics[width=0.5\textwidth, keepaspectratio]{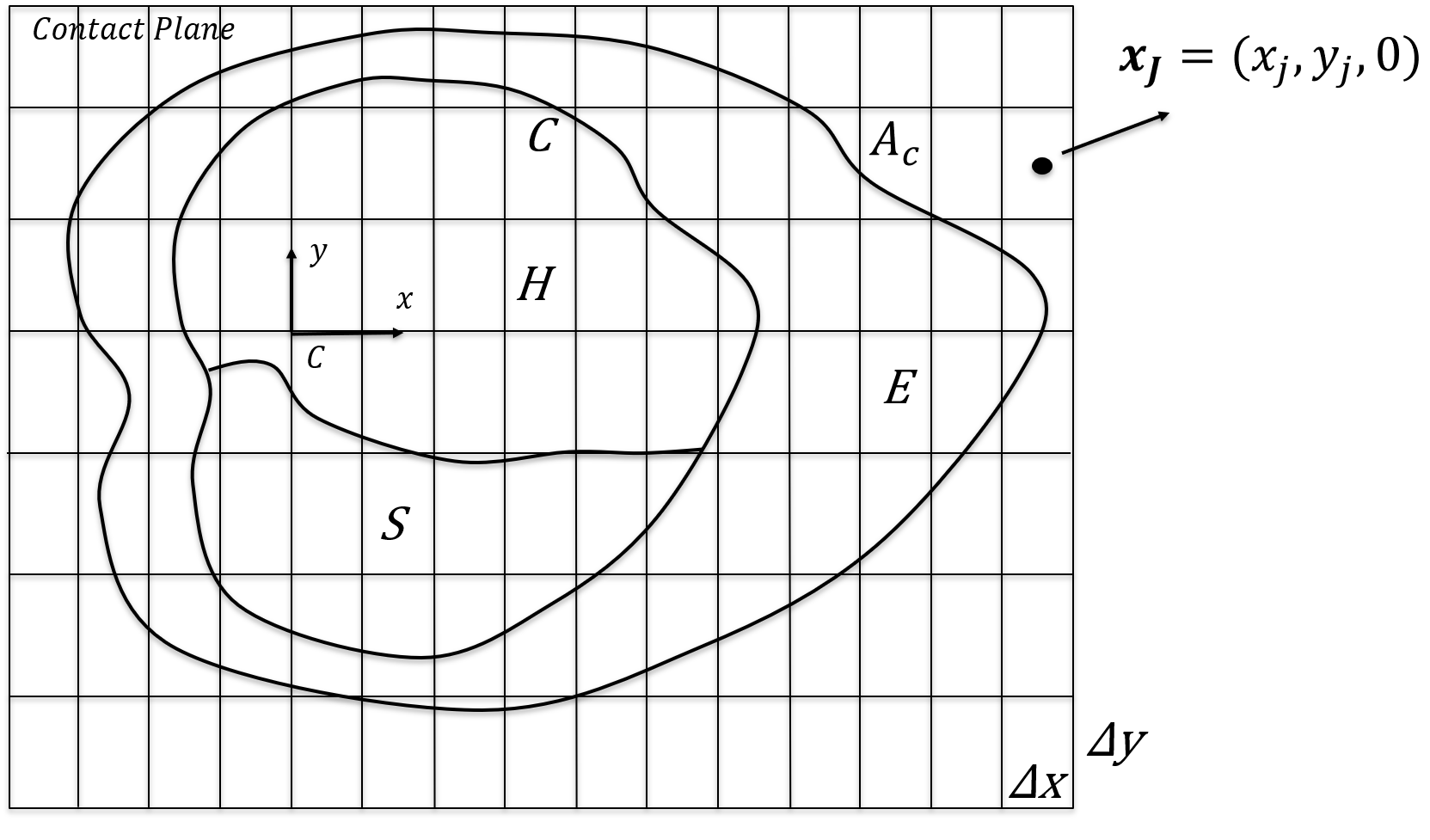}
	\caption{Local representation of the discretized contact area.}
	\label{fig:improntadiscretizzata}
\end{figure}

We give an overview of the model and algorithm used to generate our set of nonlinear systems.
Let bold letters represent vectors, the subscript $T$ denote  a  vector with components in the tangential $x$-$y$ contact place, the subscript $N$ denote  the component of a vector in the normal $z$ contact direction. 
The contact problem between two elastic bodies \cite{Kalker1, Kalker2} determines
the contact region $C$ inside the potential contact area $A_c$ (usually the interpenetration area between the wheel and rail contact surfaces), its subdivision into adhesion area $H$ and slip area $S$, and the tangential $\mathbf{p}_{T}$ and normal $p_{N}$ pressures such that the following contact conditions are satisfied:
\begin{equation}
\begin{array}{lll}
\text{ normal problem} & {\text { in contact } C: } & {e = 0, \quad p_{N} \geq 0} \\
                       & {\text { in exterior } E: } & {p_N = 0, \quad e>0} \\
                       & {C \cup E=A_{c},} & {C \cap E=\emptyset}\\
\text{ tangential problem} & {\text { in adhesion } H:} & {\| \mathbf{s_{T}}\|=0, \quad\|\mathbf{p_{T}}\| \leq g} \\
                          &{\text { in slip } S:} & {\|\mathbf{s_{T}}\| \neq 0, \quad \mathbf{p}_{T}=-g\, \mathbf{s_{T}} / \|\mathbf{s}_{T}}\| \\
&{S \cup H=C,} & {S \cap H=\emptyset} \\
\end{array}
\label{contact_problem}
\end{equation}
Above, $e$ is the deformed distance between the two bodies
and, by definition, it holds $e=0$ and  $p_N \ge 0$ in $C$.
Referring to Figure \ref{fig:improntadiscretizzata}, the region $E$  where $e>0$ is called the exterior area and  $p_N = 0$ therein.
The potential contact area is such that $A_{c}=C \cup E$.
The contact area $C$ is divided into the area of adhesion $H$ where 
the tangential component $ \mathbf{s}_{T}$ of the slip  vanishes, and the area $S$ of slip  where 
$ \mathbf{s}_{T}$ is nonzero.
The slip $ \mathbf{s}_{T}$ is the difference between the velocities of two homologous points belonging to deformed wheel and rail surfaces inside the contact area and is a 
function of the pressures $\mathbf{p}_{T}$ and $p_{N}$,  $g$ is the traction bound (Coulomb friction model \cite{Kalker1, Kalker2}).
Overall, the first three equations in (\ref{contact_problem}) model the normal contact problem 
(computation of $p_{N}$ and of the shapes of the regions $C$ and $E$), 
whereas the last three equations describe the tangential contact problem (computation of $\mathbf{p}_{T}$, of local slidings $ \mathbf{s}_{T}$ and of the shapes of the regions $H$ and $S$).

Let us consider the discretization of (\ref{contact_problem}).
 Assuming  that the contact patch is entirely contained in a plane, 
the region within which the potential contact area $A_c$ can be located is
easily discretized through a planar quadrilateral mesh, see Figure \ref{fig:improntadiscretizzata}. 
The coordinates of the center of each quadrilateral element are denoted $\mathbf{x}_{I}=\left(x_{I1}, x_{I2}, 0\right)$ 
where the capital index $I$ identifies the specific element, say $I=1, \dots, N_E$. Also, the standard indices $i=1,2,3$, will indicate the vector components. For any element $I$ and any generic vector $\mathbf{w}_I=(w_{I1},w_{I2}, w_{I3})$ associated to such mesh element, $w_{I1},w_{I2}$ are the components in the
$x$-$y$ contact plane and $w_{I3}$ is the component in the normal contact direction $z$. Namely, $\mathbf{w}_{I, T} = (w_{I1},w_{I2})$ and $w_{I3}$ are the
discrete counterparts of $\mathbf{w}_T$ and $w_N$, respectively.

The discrete values of the elastic deformation $\mathbf{u}$  
on the mesh nodes (i.e. the deformation of the elastic bodies in the contact area \cite{Kalker1, Kalker2}) are defined both at the current time instance $t$ and at the previous time instance  $t'$:
\begin{equation}
\mathbf{u}_{I} =\left(u_{I i}\right)  \hspace{0.1cm} \text { at } \hspace{0.1cm} \left(\mathbf{x}_{I}, t\right), \hspace{0.4cm}
\mathbf{u}_{I}^{\prime} =\left(u_{I i}^{\prime}\right ) \hspace{0.1cm} \text { at } \hspace{0.1cm} \left(\mathbf{x}_{I}+\mathbf{v}\left(t-t^{\prime}\right), t^{\prime}\right),
\end{equation}
where $\mathbf{v}$ is the rolling velocity (i.e. the longitudinal velocity of the wheel) and $I$ is an arbitrary mesh element). 
Analogously, for the contact pressures $\mathbf{p}$ it holds 
\begin{equation}
\mathbf{p}_{J} =\left(p_{J j}\right) \hspace{0.1cm} \text { at } \hspace{0.1cm} \left(\mathbf{x}_{J}, t\right), \hspace{0.4cm}
\mathbf{p}_{J}^{\prime} =\left(p_{J j}^{\prime}\right)  \hspace{0.1cm} \text { at } \hspace{0.1cm} \left(\mathbf{x}_{J}+\mathbf{v}\left(t-t^{\prime}\right), t^{\prime}\right),
\end{equation}
where  $J$ is an arbitrary mesh element.
According to the Boundary Element Method Theory \cite{Kalker1, Kalker2}, the discretized displacements $\mathbf{u}_{I}$ can now be written 
as a function of the discretized contact pressures $\mathbf{p}_{J}$ through the discretized version of the problem shape functions, that is
\begin{equation*}
 u_{I i} =\sum_{J=1}^{N_E} \sum_{j=1}^{3} A_{I i J j} p_{J j}, \hspace{0.2cm}  \text { with } 
A_{I i J j} :=B_{ i J j}\left(\mathbf{x}_{I}\right),
\end{equation*}
and $B_{ i J j}(\mathbf{x}_{I})$ are the discrete shape functions of the problem describing the effect of
a contact pressure $\mathbf{p}_{J}$ applied to the element $J$ on displacement $\mathbf{u}_{I}$ of the node $I$ (see \cite{Kalker1, Kalker2}). The shape function $B_{ i J j}$ usually depends on the problem geometry and the characteristics of the materials.
An analogous expression can be derived for  $u_{I i}^{\prime}$.
The elastic penetration $e$ can be calculated at each node $\mathbf{x}_I$ as
$$
e_I = h_I + \sum_J A_{I 3 J 3} p_{J 3},
$$
where $h_I$ is the discretization of the (known) undeformed distance between the two bodies, see \cite{Kalker1, Kalker2}.
Similarly, the slip $\mathbf{s}_T$ can be discretized by setting
\begin{equation}\label{sliding}
 \mathbf{s}_{I, T} = \mathbf{c}_{I, T} + (\mathbf{u}_{I, T}-\mathbf{u}_{I,T}')/(t-t'),
\end{equation}
where $\mathbf{c}_{I, T}$ is the discretization of the (given) rigid creep, that is the 
difference between the velocities of two homologous points belonging to the undeformed wheel and rail surfaces inside the contact area and thought of as rigidly connected to the bodies.  

We observe that both $\mathbf{u}$ and  $\mathbf{s}_T$ depend linearly on 
the pressures  $\mathbf{p}$ and $\mathbf{p}'$.
Therefore, the discretization of equation $e=0$ in the norm problem  (\ref{contact_problem}) yields a linear system in the discretized normal pressures $(p_{I3})$ while
the discretization of the nonlinear equation
$$
 \mathbf{p}_{T}=-g\, \mathbf{s}_{T} / \| \mathbf{s}_{T}\|,
$$
in the tangential problem  yields the nonlinear system
\begin{equation}\label{nleq_d}
 \mathbf{s}_{I, T} = - \|\mathbf{s}_{I, T}\| \mathbf{p}_{I, T} / g_I ,
\end{equation}
with $ \mathbf{p}_{I, T} = (p_{I 1}, p_{I 2})$ being the unknown\footnote{In the unlikely event  $\mathbf{s}_{I, T} = 0$, the system in nonsmooth. We regularize   (\ref{nleq_d}) replacing
the term $\sqrt{s_{I1}^2+s_{I2}^2}$ with  $\sqrt{s_{I1}^2+s_{I2}^2 + \epsilon}$, for some small positive $\epsilon$.}.
When using the Coulomb-like friction model \cite{Kalker1, Kalker2}, 
the friction limit function takes the form
$g_I = f_I p_{I3}$, where $f_I$ is a given constant friction value.

\begin{figure}[htbp]
	\centering
	  \includegraphics[width=0.5\textwidth, keepaspectratio]{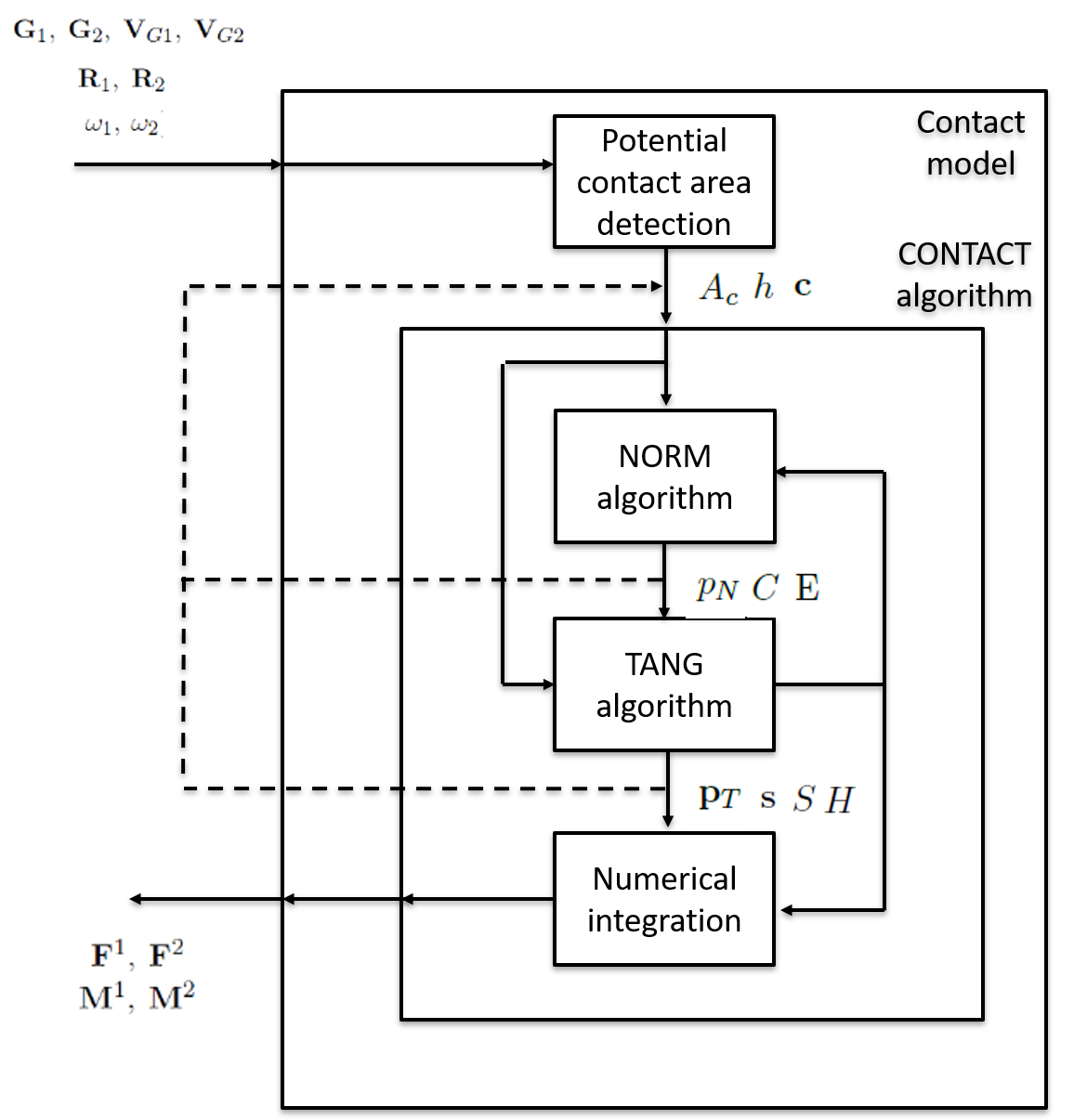}
	\caption{The architecture of the Kalker's {\small \rm CONTACT} algorithm.}
	\label{fig:archkca}
\end{figure}

The flow of Kalker's CONTACT algorithm is displayed in Figure \ref{fig:archkca} \cite{Kalker1, Kalker2}.
At each time step of  time integration, the inputs of the {\small{\rm CONTACT}} algorithm are the potential contact area $A_c$ (usually the interpenetration area between wheel and rail surfaces), 
the rigid penetration $h$ and the rigid local sliding $\mathbf{c}_T$ (inputs calculated, on turn, from the kinematic variables of the body:
position and velocities of the gravity centers $\mathbf{G}_{1}$, $\mathbf{G}_{2}$, $\mathbf{V}_{G1}$, $\mathbf{V}_{G2}$, 
rotation matrices $\mathbf{R}_{1}$, $\mathbf{R}_{2}$ and angular velocities $\mathbf{\omega}_{1}$, $\mathbf{\omega}_{2}$) \cite{Kalker1, Kalker2}. All these kinematic quantities are calculated at each time step by the ODE solver of the Simpack Rail multibody environment \cite{Simpack}.
{\small{\rm NORM}} algorithm solves the normal contact problem and returns the contact area $C$, the non-contact area $E$, the
normal contact pressures $p_N$. Then,  {\small{\rm TANG}} algorithm returns the sliding area $S$, adhesion area $H$, 
the tangential contact pressures $\mathbf{p}_T$ and local sliding $\mathbf{s}_T$. Repetitions of {\small{\rm NORM}} and {\small{\rm TANG}} 
algorithms are then performed to approximate accurately normal and tangential pressures $\mathbf{p}_T$, $p_N$. 
At the end of {\small{\rm CONTACT}} algorithm,  forces and torques exchanged by the contact bodies ($\mathbf{F}^{1}$, $\mathbf{F}^{2}$ and $\mathbf{M}^{1}$, $\mathbf{M}^{2}$) are computed by numerical integration and returned to the time integrator for proceeding in the dynamic simulation of the multibody system.

\section{Complete results}\label{sec::fullres}
In this section we collect the complete runs which gave rise  to the performance profiles in Figure \ref{Perf}. 
Results concern two velocities ($ {v}=10\, m/s$  in Tables \ref{tab:v10_rett}-\ref{tab:v10_cur}
and $ {v}=16\, m/s$ in Tables \ref{tab:v16_rett}-\ref{tab:v16_cur})
and the three different track sections (straight line in Tables \ref{tab:v10_rett} and \ref{tab:v16_rett}, cycloid in
Tables \ref{tab:v10_cy} and \ref{tab:v16_cy} and curve in
Tables \ref{tab:v10_cur} and \ref{tab:v16_cur}).
Given a sequence of nonlinear systems, we label a single system from the sequence as  Time\_Citer\_Titer   specifying
the instant time (Time), the {\small{\rm CONTACT}} iteration (Citer) and the {\small{\rm TANG}} iteration (Titer).
For each \name variant applied to a system,  we report the number of $F$-evaluations performed in case of convergence,
or, in case of failure, 
the corresponding flag.  We recall from Section \ref{numer} that a run is successful when 
$\| F_k  \| \le 10^{-6}$. 
A failure is declared either because the assigned maximum number of iterations or $F$-evaluations or backtracks is reached,
or because $\|F\|$ was not reduced for 50  consecutive  iterations.
Such occurrences are denoted as \Fit \Ffmax, \Fsigma, \Fincr, respectively.


\begin{table} 
\small
\centering
\begin{tabular}{ rrrrrrrrr }
\toprule                    
\multicolumn{ 9}{c}{$v=10\ m/s$  - straight line} \\

 System          &        \BBu & \BBd &  \BBalt& \multicolumn{ 2}{c}{\small \rm ABB} & \multicolumn{ 2}{c}{\small \rm ABBm} &  \dabbmino  \\

           &            & &           & $\tau=0.1$ & $\tau=0.8$ & $\tau=0.1$ & $\tau=0.8$ &         \\
\midrule  
101\_1\_2&	69&	59&	74&	75&	59&	71&	57&	69\\
101\_2\_2&	382&	148&	248&	295&	205&	174&	198&	220\\
103\_1\_2&	37&	31&	35&	37&	30&	37&	31&	34\\
103\_2\_2&	37&	31&	35&	37&	30&	37&	31&	34\\
104\_1\_2&	36&	36&	37&	36&	38&	36&	39&	38\\
104\_2\_2&	36&	36&	37&	36&	38&	36&	39&	38\\
105\_1\_2&	39&	38&	39&	39&	38&	39&	39&	39\\
105\_1\_3&	77&	69&	82&	79&	70&	82&	67&	74\\
105\_2\_2&	40&	37&	39&	40&	38&	40&	39&	39\\
105\_2\_3&	74&	73&	86&	75&	70&	75&	67&	76

\\

\bottomrule
\end{tabular}

\caption{Number of function evaluations performed by \name variants in the solution of nonlinear systems arising from time 100 to time 105 and corresponding to a straight line with velocity $10\ m/s$. In the first column we indicate the time step, the {\small{\rm CONTACT}} and the {\small{\rm TANG}} iteration.}\label{tab:v10_rett}
\end{table}

\begin{sidewaystable} 
\centering
\footnotesize

\begin{tabular}{rrrrrrrrr rrrrrrrrr}
\toprule
           &            \multicolumn{ 17}{c}{velocity $10\ m/s$  - cycloid} \\

 System          &        \BBu & \BBd & \BBalt & \multicolumn{ 2}{c}{\small \rm ABB} & \multicolumn{ 2}{c}{\small \rm ABBm} &  \dabbmino &   System          &        \BBu & \BBd &  \BBalt & \multicolumn{ 2}{c}{\small \rm ABB} & \multicolumn{ 2}{c}{\small \rm ABBm} &  \dabbmino \\

           &            &       &     & $\tau=0.1$ & $\tau=0.8$ & $\tau=0.1$ & $\tau=0.8$ &                       &    &        &      &      & $\tau=0.1$ & $\tau=0.8$ & $\tau=0.1$ & $\tau=0.8$          &   \\
\midrule
    300\_1\_2 & 178   & 128   & 137   & 145   & 149   & 174   & 133   & 163   & 303\_2\_2 & \Ffmax & \Fincr & 2196  & \Fincr & \Fincr & 1111  & 763   & 887 \\
    300\_1\_3 & 513   & 304   & 257   & 296   & 252   & 271   & 230   & 298   & 303\_2\_3 & \Ffmax & 1062  & 7400  & 1486  & 1413  & 911   & 722   & 798 \\
    300\_1\_4 & 569   & 402   & 290   & 464   & 350   & 460   & 278   & 299   & 303\_2\_4 & \Ffmax & 1713  & 10229 & 1780  & 1400  & \Fincr & 889   & 1054 \\
    300\_2\_2 & 343   & 203   & 266   & 229   & 194   & 209   & 168   & 204   & 303\_2\_5 & \Ffmax & 1424  & 23393 & 2053  & 1776  & 1201  & 1046  & 1358 \\
    300\_2\_3 & 16421 & 388   & 398   & 406   & 686   & 410   & 330   & 408   & 303\_3\_2 & \Ffmax & 926   & 6424  & 1352  & 806   & 896   & 814   & 821 \\
    300\_3\_2 & 357   & 223   & 248   & 257   & 205   & 225   & 187   & 232   & 303\_3\_3 & \Ffmax & 1318  & 6285  & 1508  & 886   & 1074  & 981   & 896 \\
    300\_3\_3 & 1650  & 385   & 368   & 432   & 530   & 462   & 339   & 499   & 303\_3\_4 & \Ffmax & 1279  & 14647 & 2295  & 1501  & 1244  & 959   & 1012 \\
    301\_1\_2 & 415   & 281   & 247   & 326   & 325   & 264   & 243   & 248   & 303\_3\_5 & \Ffmax & \Fincr & 17619 & 2353  & \Fincr & 1484  & 1311  & 1193 \\
    301\_1\_3 & 503   & 319   & 351   & 342   & 480   & 280   & 286   & 329   & 304\_1\_2 & 39075 & 962   & 815   & 643   & 504   & 714   & 447   & 491 \\
    301\_1\_4 & 582   & 442   & 281   & 380   & 376   & 344   & 291   & 305   & 304\_1\_3 & \Ffmax & 711   & 2891  & 860   & 1242  & 710   & 607   & 562 \\
    301\_2\_2 & 1127  & 286   & 298   & 271   & 430   & 310   & 284   & 297   & 304\_1\_4 & \Ffmax & 1524  & 3611  & 966   & 1423  & 785   & 515   & 752 \\
    301\_2\_3 & 630   & 414   & 367   & 388   & 430   & 322   & 313   & 337   & 304\_2\_2 & 725   & 366   & 381   & 393   & 416   & 300   & 311   & 317 \\
    301\_2\_4 & 758   & 345   & 372   & 408   & 355   & 363   & 319   & 386   & 304\_2\_3 & 65775 & 558   & 648   & 753   & 734   & 577   & 453   & 548 \\
    301\_3\_2 & 918   & 357   & 299   & 315   & 350   & 294   & 288   & 326   & 304\_2\_4 & 56953 & 709   & 1870  & 638   & 920   & 562   & 475   & 523 \\
    301\_3\_3 & 750   & 400   & 320   & 473   & 423   & 350   & 305   & 313   & 304\_3\_2 & 415   & 421   & 370   & 470   & 431   & 357   & 339   & 325 \\
    301\_3\_4 & 440   & 363   & 302   & 352   & 434   & 310   & 301   & 393   & 304\_3\_3 & 47176 & 533   & 2376  & 616   & 627   & 518   & 411   & 612 \\
    302\_1\_2 & \Ffmax & 743   & 3727  & 993   & 1022  & 558   & 457   & 495   & 304\_3\_4 & 86605 & 696   & 1180  & 709   & 603   & 557   & 468   & 488 \\
    302\_1\_3 & \Ffmax & 844   & 4067  & 1183  & 972   & 1068  & 670   & 678   & 305\_1\_2 & 796   & 270   & 311   & 302   & 323   & 329   & 242   & 364 \\
    302\_1\_4 & \Ffmax & 3546  & 25810 & 6171  & 2529  & 1735  & 1267  & 1342  & 305\_1\_3 & 339   & 293   & 270   & 271   & 294   & 288   & 243   & 310 \\
    302\_2\_2 & 634   & 444   & 417   & 552   & 539   & 431   & 332   & 376   & 305\_1\_4 & 430   & 342   & 301   & 354   & 335   & 307   & 230   & 309 \\
    302\_2\_3 & 27285 & 610   & 508   & 890   & 544   & 502   & 398   & 548   & 305\_2\_2 & \Ffmax & \Fincr & 2434  & 1401  & 800   & \Fincr & 1282  & 1208 \\
    302\_2\_4 & \Ffmax & \Fincr & 7325  & 1359  & 1951  & 927   & 853   & 693   & 305\_2\_3 & \Ffmax & 1110  & 2222  & 1713  & 1030  & 950   & 717   & 684 \\
    302\_3\_2 & 743   & 426   & 373   & 455   & 438   & 402   & 332   & 361   & 305\_2\_4 & \Ffmax & \Fincr & 842   & 1527  & 846   & 748   & 768   & 648 \\
    302\_3\_3 & 39825 & 739   & 502   & 869   & 616   & 459   & 401   & 463   & 305\_2\_5 & \Ffmax & \Fincr & 3329  & 1516  & 850   & 1332  & 573   & 597 \\
    302\_3\_4 & \Ffmax & 2245  & 7598  & 1141  & 938   & 1005  & 660   & 702   & 305\_3\_2 & \Ffmax & 980   & 6755  & 1524  & \Fincr & 920   & 1036  & 1518 \\
    303\_1\_2 & 22687 & 554   & 679   & 502   & \Fincr & 609   & 405   & 460   & 305\_3\_3 & \Ffmax & \Fincr & 5805  & 1829  & 756   & 694   & 634   & 579 \\
    303\_1\_3 & 33798 & 468   & 684   & 571   & 578   & 461   & 411   & 562   & 305\_3\_4 & \Ffmax & 871   & 2502  & 1363  & 997   & 857   & 716   & 648 \\
    303\_1\_4 & \Ffmax & 965   & 1163  & 734   & 669   & 653   & 524   & 613   & 305\_3\_5 & \Ffmax & \Fincr & 1786  & 1286  & 843   & 929   & 702   & 663 \\
    \bottomrule
    \end{tabular}

\caption{Results for each system of the sequences generated in the cycloid section of the train track with velocity ${v}=10\ m/s$.}\label{tab:v10_cy}
\end{sidewaystable} 

\begin{sidewaystable} 
\centering
\footnotesize

\begin{tabular}{rrrrrrrrr rrrrrrrrr}
\toprule           &                   \multicolumn{ 17}{c}{velocity $10\ m/s$  - curve} \\

 System          &        \BBu & \BBd &  \BBalt & \multicolumn{ 2}{c}{\small \rm ABB} & \multicolumn{ 2}{c}{\small \rm ABBm} &  \dabbmino &  System          &        \BBu & \BBd &  \BBalt & \multicolumn{ 2}{c}{\small \rm ABB} & \multicolumn{ 2}{c}{\small \rm ABBm} &  \dabbmino \\

           &            &            & & $\tau=0.1$ & $\tau=0.8$ & $\tau=0.1$ & $\tau=0.8$ &           &            &     &       &            & $\tau=0.1$ & $\tau=0.8$ & $\tau=0.1$ & $\tau=0.8$ &           \\
\midrule

    450\_1\_2 & 386   & 210   & 246   & 251   & 293   & 293   & 211   & 284   & 453\_1\_3 & 402   & 319   & 457 &427   & 405   & 409   & 255   & 316 \\
    450\_1\_3 & 623   & 204   & 303   & 285   & 281   & 268   & 1580  & 1627  & 453\_1\_4 & \Ffmax & \Fincr & 2705 &656 & 1285  & 996   & 611   & 544 \\
    450\_2\_2 & 29520 & 492   & 457   & 475   & 416   & 458   & 320   & 471   & 453\_2\_2 & 536   & 356   & 379 &593  & 409   & 362   & 329   & 355 \\
    450\_2\_3 & 12031 & 428   & 433   & 412   & 458   & 415   & 309   & 387   & 453\_2\_3 & \Ffmax & 739   & 872 &1030   & 557   & 726   & \Fincr & 560 \\
    450\_3\_2 & 13652 & 560   & 403   & 562   & 416   & 463   & 379   & 382   & 453\_2\_4 & \Ffmax & 1772  & \Fincr &\Fincr & 2018  & 1579  & 1535  & \Fincr \\
    450\_3\_3 & 11509 & 464   & 448   & 518   & 493   & 475   & 393   & 391   & 453\_3\_2 & 566   & 351   & 355 &548   & 392   & 367   & 337   & 398 \\
    451\_1\_2 & 681   & 437   & 382   & 520   & 570   & 519   & 340   & 397   & 453\_3\_3 & \Ffmax & 558   & 598 &796   & 617   & 612   & 536   & 568 \\
    451\_1\_3 & \Ffmax & 1218  & 4314  & 999   & 1564  & 868   & 613   & 1501  & 453\_3\_4 & \Ffmax & \Fincr & \Fsigma &2308 & \Fincr & 1487  & 1187  & 1667 \\
    451\_1\_4 & \Ffmax & 3805  & 18920 & 1790  & \Fincr & 1305  & 1083  & 1334  & 454\_1\_2 & 147   & 153   & 165 &139   & 153   & 137   & 138   & 150 \\
    451\_2\_2 & 324   & 274   & 329   & 264   & 264   & 263   & 210   & 250   & 454\_1\_3 & 207   & 175   & 206  &229 & 192   & 194   & 154   & 175 \\
    451\_2\_3 & \Ffmax & 1652  & 1046  & 859   & 1304  & 691   & 520   & 595   & 454\_1\_4 & 2367  & 276   & 293 &286  & 332   & 283   & 252   & 314 \\
    451\_2\_4 & \Ffmax & 1573  & \Fincr & 1260  & \Fincr & 1232  & \Fincr & 941   & 454\_1\_5 & 861   & 351   & 250 &269  & 332   & 291   & 231   & 301 \\
    451\_3\_2 & 381   & 253   & 240   & 301   & 243   & 285   & 209   & 270   & 454\_2\_2 & 237   & 172   & 209 &194  & 191   & 202   & 153   & 207 \\
    451\_3\_3 & \Ffmax & 3141  & 4232  & 660   & 801   & 640   & 606   & 635   & 454\_2\_3 & 413   & 279   & 211 &288   & 315   & 240   & 254   & 280 \\
    451\_3\_4 & \Ffmax & \Fincr & \Fincr & \Fincr & \Fincr & 1042  & 936   & 888   & 454\_2\_4 & 901   & 363   & 209 &256   & 307   & 262   & 227   & 261 \\
    451\_4\_2 & 358   & 296   & 321   & 279   & 295   & 268   & 213   & 263   & 454\_3\_2 & 259   & 204   & 204  &183 & 198   & 183   & 157   & 183 \\
    451\_4\_3 & \Ffmax & 2108  & 901   & 688   & 729   & 676   & 597   & 639   & 454\_3\_3 & 469   & 317   & 329 &273   & 290   & 244   & 251   & 265 \\
    451\_4\_4 & \Ffmax & \Fincr & 12872 & 1797  & \Fincr & 1093  & 905   & 821   & 454\_3\_4 & 450   & 302   & 231 &277  & 297   & 254   & 229   & 270 \\
    452\_1\_2 & 66785 & 638   & 638   & 548   & 743   & 585   & 545   & 522   & 455\_1\_2 & 147   & 137   & 145 &144  & 126   & 145   & 127   & 136 \\
    452\_1\_3 & 71198 & 701   & 725   & 535   & 789   & 489   & 552   & 508   & 455\_1\_3 & 212   & 184   & 203 &219  & 166   & 226   & 166   & 196 \\
    452\_1\_4 & 45680 & 803   & 521   & 617   & 594   & 584   & 470   & 520   & 455\_1\_4 & 482   & 272   & 256 &291  & 278   & 251   & 237   & 246 \\
    452\_2\_2 & 498   & 557   & 887   & 514   & 539   & 417   & 301   & 467   & 455\_2\_2 & 497   & 372   & 250 &496   & 288   & 256   & 270   & 284 \\
    452\_2\_3 & 37679 & 608   & 714   & 474   & 672   & 456   & 425   & 454   & 455\_2\_3 & 563   & 393   & 473 &641   & 340   & 436   & 357   & 348 \\
    452\_2\_4 & 40269 & 718   & 797   & 565   & 790   & 484   & 379   & 501   & 455\_2\_4 & \Ffmax & 840   & 5928 &1544 & 929   & 1131  & 618   & 632 \\
    452\_3\_2 & 31230 & 433   & 451   & 438   & 517   & 345   & 405   & 354   & 455\_3\_2 & 341   & 270   & 268 &391  & 392   & 302   & 238   & 282 \\
    452\_3\_3 & 41623 & 581   & 634   & 575   & 726   & 509   & 400   & 451   & 455\_3\_3 & 603   & 432   & 405  &592 & 415   & 363   & 346   & 353 \\
    452\_3\_4 & 5592  & 477   & 658   & 572   & 570   & 457   & 407   & 470   & 455\_3\_4 & \Ffmax & 792   & 7505 &1586 & 855   & 914   & 663   & 744 \\
    453\_1\_2 & 288   & 200   & 257   & 227   & 210   & 279   & 190   & 210   &       &       &       &       &       &       &       &  \\
    \bottomrule
    \end{tabular}%

\caption{Results for each system of the sequences generated in the curve segment of the train path with velocity ${v}=10\ m/s$.}\label{tab:v10_cur}
\end{sidewaystable}

\begin{table} 
\small 
\centering
\begin{tabular}{rrrrrrrrr}

\toprule   

\multicolumn{ 9}{c}{velocity $16\ m/s$  - straight line} \\

   System          &         \BBu & \BBd &   \BBalt & \multicolumn{ 2}{c}{\small \rm ABB} & \multicolumn{ 2}{c}{\small \rm ABBm} &  \dabbmino \\

           &            &          &  & $\tau=0.1$ & $\tau=0.8$ & $\tau=0.1$ & $\tau=0.8$ &           \\
\midrule
    50\_1\_2&	60&	45&	53&	52&	47&	52&	46&	49\\
50\_2\_2&	53&	44&	51&	54&	48&	54&	48&	53\\
50\_3\_2&	53&	44&	51&	48&	48&	48&	48&	53\\
52\_2\_2&	75&	78&	53&	76&	75&	101&	61&	91\\
52\_3\_2&	89&	78&	53&	76&	88&	112&	61&	91\\
55\_1\_2&	65&	66&	66&	83&	66&	80&	62&	72\\
55\_2\_2&	69&	79&	60&	76&	61&	73&	67&	71\\
55\_3\_2&	69&	79&	60&	80&	61&	73&	67&	71

    \\
\bottomrule
\end{tabular}

\caption{Number of function evaluations performed by \name variants in the solution of nonlinear systems arising from time 50 to time 55 and corresponding to a straight line with velocity $16\ m/s$. In the first column we indicate the time step, the 
{\small{\rm CONTACT}} and the {\small{\rm TANG}} iteration.}\label{tab:v16_rett}
\end{table}

\begin{sidewaystable} 
\centering
\footnotesize
\begin{tabular}{rrrrrrrrr rrrrrrrrr}
\toprule
           &                        \multicolumn{ 17}{c}{velocity $16\ m/s $  - cycloid} \\

 System          &        \BBu & \BBd & \BBalt & \multicolumn{ 2}{c}{\small \rm ABB} & \multicolumn{ 2}{c}{\small \rm ABBm} &  \dabbmino &   System          &        \BBu & \BBd &  \BBalt  & \multicolumn{ 2}{c}{\small \rm ABB} & \multicolumn{ 2}{c}{\small \rm ABBm} &  \dabbmino  \\
           
   &            &      &      & $\tau=0.1$ & $\tau=0.8$ & $\tau=0.1$ & $\tau=0.8$             &            &   &         &      &      & $\tau=0.1$ & $\tau=0.8$ & $\tau=0.1$ & $\tau=0.8$ &           \\   
   \midrule
  
    150\_1\_2 & 985   & 297   & 330   & 366   & 357   & 351   & 278   & 343   & 153\_1\_3 & \Ffmax & 1173  & 1181  & 1162  & 1179  & 735   & 568   & 596 \\
    150\_1\_3 & 26886 & 569   & 512   & 612   & 555   & 487   & 419   & 437   & 153\_1\_4 & \Ffmax & 991   & 3881  & 1003  & 1590  & 1044  & 635   & 771 \\
    150\_1\_4 & \Ffmax & 967   & 3163  & 653   & \Fincr & 550   & 604   & 617   & 153\_2\_2 & 21846 & 475   & 603   & 688   & 532   & 578   & 396   & 446 \\
    150\_1\_5 & \Ffmax & \Fincr & 810   & 647   & 1549  & 614   & 510   & 710   & 153\_2\_3 & \Ffmax & 1149  & 3920  & 1316  & 1506  & 843   & 621   & 704 \\
    150\_2\_2 & 476   & 228   & 307   & 295   & 302   & 277   & 216   & 301   & 153\_2\_4 & \Ffmax & 1445  & 5035  & 1262  & 1272  & 1215  & 602   & 784 \\
    150\_2\_3 & 627   & 584   & 404   & 437   & 485   & 377   & 344   & 443   & 153\_2\_5 & \Ffmax & 772   & 4023  & 926   & 1576  & 1188  & 764   & 725 \\
    150\_2\_4 & 52373 & 585   & 479   & 494   & 730   & 438   & 391   & 435   & 153\_3\_2 & 1873  & 628   & 754   & 674   & 585   & 489   & 429   & 471 \\
    150\_3\_2 & \Ffmax & 1304  & \Fincr & \Fincr & 1777  & 2707  & 1237  & 911   & 153\_3\_3 & \Ffmax & 770   & 4768  & 1187  & 1882  & 941   & 699   & 860 \\
    150\_3\_3 & \Ffmax & 2498  & \Fincr & \Fincr & \Fincr & 2300  & 1973  & 1737  & 153\_3\_4 & \Ffmax & 1568  & 4872  & 923   & 1161  & 1173  & 678   & 709 \\
    150\_3\_4 & \Ffmax & 6214  & \Fincr & \Fincr & \Fincr & 3097  & 2576  & \Fincr & 153\_3\_5 & \Ffmax & 1226  & 5474  & 1145  & 1118  & 730   & 688   & 730 \\
    151\_1\_2 & \Ffmax & \Fincr & 5095  & 841   & 905   & 664   & 605   & 689   & 154\_1\_2 & 66851 & 776   & 3124  & 727   & 1033  & 585   & 534   & 527 \\
    151\_1\_3 & \Ffmax & 1114  & 5312  & 1421  & 1144  & 810   & 616   & 829   & 154\_1\_3 & 1031  & 386   & 513   & 467   & 681   & 433   & 310   & 346 \\
    151\_1\_4 & \Ffmax & 1454  & 8154  & 1630  & 3755  & 1125  & 1139  & 1046  & 154\_1\_4 & 18703 & 533   & 421   & 539   & 518   & 434   & 404   & 447 \\
    151\_1\_5 & \Ffmax & 3590  & 13111 & 2610  & 1435  & 1231  & 864   & 1043  & 154\_2\_2 & 947   & 319   & 312   & 420   & 357   & 341   & 294   & 356 \\
    151\_2\_2 & \Ffmax & 1337  & 12656 & 1333  & 3092  & 973   & 864   & 856   & 154\_2\_3 & 255   & 193   & 220   & 216   & 241   & 238   & 201   & 246 \\
    151\_2\_3 & \Ffmax & 3776  & 9599  & 1983  & 2198  & 1077  & 949   & 961   & 154\_2\_4 & 348   & 266   & 255   & 255   & 258   & 250   & 228   & 276 \\
    151\_2\_4 & \Ffmax & 3013  & 9073  & 1867  & 3551  & 1409  & 870   & 974   & 154\_3\_2 & 569   & 403   & 288   & 336   & 394   & 302   & 277   & 354 \\
    151\_2\_5 & \Ffmax & 5005  & 18543 & 1831  & 3662  & 1635  & 1270  & 1345  & 154\_3\_3 & 248   & 218   & 249   & 253   & 276   & 217   & 206   & 233 \\
    151\_3\_2 & \Ffmax & \Fincr & 7743  & \Fincr & 3893  & \Fincr & 939   & 803   & 154\_3\_4 & 346   & 318   & 278   & 281   & 271   & 267   & 239   & 250 \\
    151\_3\_3 & \Ffmax & 2293  & 9494  & 1383  & 1689  & 1080  & 809   & 982   & 155\_1\_2 & \Ffmax & 1161  & 5470  & 1151  & 987   & 824   & 718   & 859 \\
    151\_3\_4 & \Ffmax & 1235  & 7622  & 1416  & 1884  & 1075  & 856   & 941   & 155\_1\_3 & \Ffmax & \Fincr & 31313 & 4192  & 4270  & 1758  & 1401  & 1193 \\
    151\_3\_5 & \Ffmax & 4085  & 24983 & 1853  & \Fincr & 1509  & 1147  & 1330  & 155\_1\_4 & \Ffmax & 5839  & 19894 & \Fincr & 4182  & 1621  & 1729  & 1380 \\
    152\_1\_2 & 68856 & 822   & 1395  & 742   & 661   & 680   & 473   & 575   & 155\_1\_5 & \Ffmax & \Fincr & \Fincr & \Fincr & \Fincr & 1624  & 1351  & 1339 \\
    152\_1\_3 & \Ffmax & 682   & 4009  & 1153  & 1085  & 859   & 648   & 669   & 155\_2\_2 & \Ffmax & 1211  & 3754  & 1267  & 1275  & 764   & 651   & 635 \\
    152\_1\_4 & \Ffmax & 725   & 2905  & 986   & 1423  & 799   & 646   & 720   & 155\_2\_3 & \Ffmax & \Fincr & \Fincr & 2536  & \Fincr & 1658  & 1328  & 1273 \\
    152\_2\_2 & 21104 & 604   & 641   & 407   & 681   & 543   & 347   & 399   & 155\_2\_4 & \Ffmax & 1623  & 24770 & 3690  & \Fincr & 1626  & 1461  & 1427 \\
    152\_2\_3 & 80349 & 701   & 1082  & 636   & 845   & 632   & 476   & 610   & 155\_2\_5 & \Ffmax & \Fincr & \Fsigma & \Fincr & \Fincr & 1683  & 1715  & 1559 \\
    152\_2\_4 & \Ffmax & 1748  & 3725  & 1395  & 1034  & 873   & 590   & 849   & 155\_3\_2 & \Ffmax & 877   & 6004  & 990   & 882   & 795   & 567   & 818 \\
    152\_3\_2 & 20711 & 567   & 601   & 382   & 664   & 453   & 358   & 420   & 155\_3\_3 & \Ffmax & \Fincr & 23302 & 1784  & \Fincr & \Fincr & 1539  & 1238 \\
    152\_3\_3 & 75894 & 966   & 1098  & 522   & 898   & 639   & 535   & 627   & 155\_3\_4 & \Ffmax & 2895  & 32130 & 1953  & \Fincr & 1539  & 1739  & 1315 \\
    152\_3\_4 & \Ffmax & 1146  & 4114  & 848   & 1152  & 744   & 558   & 734   & 155\_3\_5 & \Ffmax & \Fincr & \Fincr & 6554  & \Fincr & \Fincr & \Fincr & \Fincr \\
    153\_1\_2 & 1281  & 408   & 589   & 512   & 495   & 472   & 400   & 397   &       &       &       &       &       &       &       &       &  \\
    \bottomrule
    \end{tabular}

\caption{Results for each system of the sequences generated in the cycloid section
of the train track with velocity ${v}=16\ m/s$.}
\label{tab:v16_cy}
\end{sidewaystable}

\begin{sidewaystable} 
\centering
\footnotesize
\begin{tabular}{rrrrrrrrr rrrrrrrrr}
\toprule
           &              \multicolumn{ 17}{c}{velocity $16\ m/s $ - curve} \\

 System          &        \BBu & \BBd & \BBalt & \multicolumn{ 2}{c}{\small \rm ABB} & \multicolumn{ 2}{c}{\small \rm ABBm} &  \dabbmino &    System          &        \BBu & \BBd &  \BBalt & \multicolumn{ 2}{c}{\small \rm ABB} & \multicolumn{ 2}{c}{\small \rm ABBm} &  \dabbmino  \\
           
   &            &       &     & $\tau=0.1$ & $\tau=0.8$ & $\tau=0.1$ & $\tau=0.8$    &         &            &            &         &   & $\tau=0.1$ & $\tau=0.8$ & $\tau=0.1$ & $\tau=0.8$      &       \\   
   \midrule
   
    350\_1\_2 & 424   & 320   & 308   & 359   & 366   & 297   & 284   & 286   & 352\_4\_5 & \Ffmax & 1132  & 7322  & 1252  & \Fincr & 921   & \Fincr & 724 \\
    350\_1\_3 & \Ffmax & 825   & 5650  & 826   & 905   & 771   & 540   & 687   & 353\_1\_2 & 468   & 357   & 398   & 482   & 342   & 352   & 307   & 357 \\
    350\_2\_2 & 308   & 208   & 220   & 244   & 261   & 243   & 197   & 247   & 353\_1\_3 & 887   & 640   & 588   & 557   & 441   & 508   & 446   & 456 \\
    350\_2\_3 & \Ffmax & 1322  & 3384  & 572   & \Fincr & 501   & 433   & 497   & 353\_1\_4 & \Ffmax & 695   & 4525  & 905   & 1369  & 781   & 625   & 656 \\
    350\_2\_4 & \Ffmax & \Fincr & 6845  & 1204  & 1523  & 746   & 790   & 718   & 353\_1\_5 & \Ffmax & 877   & 4670  & 793   & 1551  & 782   & 682   & 764 \\
    350\_3\_2 & 311   & 221   & 277   & 264   & 234   & 214   & 188   & 213   & 353\_2\_2 & 589   & 357   & 365   & 461   & 398   & 426   & 370   & 386 \\
    350\_3\_3 & 76754 & \Fincr & 885   & 639   & 666   & 491   & 416   & 481   & 353\_2\_3 & 47619 & 755   & 572   & 913   & 812   & 529   & 459   & 528 \\
    350\_3\_4 & \Ffmax & \Fincr & 6032  & 675   & \Fincr & 1141  & 761   & 647   & 353\_2\_4 & \Ffmax & 1143  & 3476  & \Fincr & 857   & 798   & 642   & 687 \\
    350\_4\_2 & 271   & 207   & 233   & 229   & 226   & 220   & 201   & 218   & 353\_2\_5 & \Ffmax & 1984  & 8598  & 1370  & 1700  & \Fincr & 867   & 1111 \\
    350\_4\_3 & 91233 & 764   & 3110  & 633   & 829   & 536   & 432   & 526   & 353\_3\_2 & 711   & 381   & 394   & 481   & 380   & 408   & 368   & 361 \\
    350\_4\_4 & \Ffmax & 1593  & 6301  & 722   & \Fincr & 637   & \Fincr & 751   & 353\_3\_3 & 65122 & 672   & 600   & 710   & 996   & 604   & 511   & 457 \\
    351\_1\_2 & \Ffmax & 1241  & 1625  & 920   & 913   & 772   & 597   & 538   & 353\_3\_4 & \Ffmax & 837   & 1623  & 815   & 1111  & 759   & 588   & 633 \\
    351\_1\_3 & \Ffmax & 1596  & 11134 & 1807  & \Fincr & 1374  & 1199  & 1090  & 353\_3\_5 & \Ffmax & 1250  & 6524  & 1233  & 1350  & 1110  & 915   & 855 \\
    351\_1\_4 & \Ffmax & 2272  & 20207 & 1862  & \Fincr & 1555  & 1217  & 1240  & 353\_4\_2 & 575   & 448   & 505   & 425   & 360   & 350   & 341   & 372 \\
    351\_2\_2 & \Ffmax & 1088  & \Fincr & \Fincr & 1207  & 1385  & 959   & 1050  & 353\_4\_3 & 57903 & 732   & 725   & 644   & 469   & 517   & 492   & 533 \\
    351\_2\_3 & \Ffmax & 2428  & \Fincr & \Fincr & \Fincr & 2185  & 1567  & 1825  & 353\_4\_4 & \Ffmax & 1030  & 932   & 873   & 1055  & 679   & 630   & 669 \\
    351\_2\_4 & \Ffmax & 5683  & \Fincr & \Fincr & \Fincr & 2421  & 2064  & 1636  & 353\_4\_5 & \Ffmax & \Fincr & 8112  & 1276  & 1502  & 980   & 904   & 967 \\
    351\_2\_5 & \Ffmax & \Fincr & \Fincr & \Fincr & \Fincr & 3192  & 2052  & 2770  & 354\_1\_2 & 313   & 229   & 219   & 320   & 261   & 265   & 187   & 253 \\
    351\_3\_2 & \Ffmax & 1261  & 12388 & 3742  & 1566  & 992   & 1166  & 876   & 354\_1\_3 & 502   & 323   & 369   & 398   & 337   & 318   & 267   & 342 \\
    351\_3\_3 & \Ffmax & 2029  & \Fincr & \Fincr & \Fincr & \Fincr & \Fincr & 1704  & 354\_1\_4 & 87446 & 710   & 4042  & 610   & 716   & 579   & 536   & 673 \\
    351\_3\_4 & \Ffmax & 2397  & \Fincr & \Fincr & 4270  & 2105  & 2074  & 1630  & 354\_2\_2 & 445   & 321   & 348   & 373   & 292   & 289   & 230   & 296 \\
    351\_3\_5 & \Ffmax & \Fincr & \Fincr & \Fincr & \Fincr & 2833  & \Fincr & 2635  & 354\_2\_3 & 1771  & 462   & 359   & 434   & 473   & 355   & 345   & 372 \\
    351\_4\_2 & \Ffmax & 1285  & \Fincr & 4846  & 1378  & 1262  & 1313  & 1028  & 354\_2\_4 & \Ffmax & 1054  & 4522  & 1052  & 1159  & 757   & 649   & 701 \\
    351\_4\_3 & \Ffmax & 1778  & \Fincr & \Fincr & 2581  & 2073  & 2144  & 1764  & 354\_3\_2 & 451   & 315   & 295   & 324   & 275   & 259   & 265   & 316 \\
    351\_4\_4 & \Ffmax & \Fincr & \Fincr & \Fincr & \Fincr & 2848  & 1794  & 1763  & 354\_3\_3 & 789   & 382   & 392   & 508   & 521   & 409   & 408   & 409 \\
    351\_4\_5 & \Ffmax & \Fincr & \Fincr & \Fincr & \Fincr & \Fincr & 3340  & \Fincr & 354\_3\_4 & \Ffmax & 913   & 3478  & 786   & 921   & 845   & 607   & 665 \\
    352\_1\_2 & \Ffmax & 1794  & \Fsigma & 5760  & 1636  & 1619  & 1933  & 1728  & 354\_4\_2 & 405   & 323   & 289   & 350   & 308   & 317   & 256   & 295 \\
    352\_1\_3 & \Ffmax & 3141  & \Fsigma & 3787  & 2872  & 1686  & 1495  & 1524  & 354\_4\_3 & 1776  & 497   & 363   & 452   & 338   & 399   & 333   & 370 \\
    352\_1\_4 & \Ffmax & \Fincr & \Fincr & \Fincr & \Fincr & 2334  & 1657  & 1721  & 354\_4\_4 & \Ffmax & 991   & 4561  & 830   & 1141  & 704   & 553   & 634 \\
    352\_1\_5 & \Ffmax & \Fincr & \Fincr & \Fincr & \Fincr & 2318  & 2846  & 1623  & 355\_1\_2 & 638   & 226   & 262   & 264   & 292   & 268   & 258   & 266 \\
    352\_2\_2 & 72375 & 676   & 1359  & 708   & 586   & 643   & 459   & 501   & 355\_1\_3 & 527   & 339   & 509   & 348   & 348   & 348   & 286   & 331 \\
    352\_2\_3 & 74955 & 801   & 878   & 794   & 718   & 857   & 481   & 519   & 355\_1\_4 & 35134 & 489   & 1201  & 464   & 525   & 477   & 382   & 408 \\
    352\_2\_4 & \Ffmax & 866   & 5116  & 1209  & 1071  & 837   & 648   & 746   & 355\_2\_2 & 346   & 222   & 252   & 246   & 243   & 221   & 194   & 242 \\
    352\_2\_5 & \Ffmax & \Fincr & 12683 & 1209  & \Fincr & 921   & 803   & 909   & 355\_2\_3 & 2303  & 480   & 396   & 402   & 357   & 313   & 261   & 358 \\
    352\_3\_2 & 59157 & 701   & 1249  & 712   & 652   & 687   & 420   & 589   & 355\_2\_4 & 41075 & 671   & 542   & 511   & 401   & 376   & 355   & 433 \\
    352\_3\_3 & 87628 & 1116  & 682   & 804   & 611   & 639   & 517   & 517   & 355\_3\_2 & 336   & 289   & 249   & 264   & 282   & 194   & 232   & 241 \\
    352\_3\_4 & \Ffmax & 808   & 6379  & 845   & 830   & 726   & 782   & 685   & 355\_3\_4 & 639   & 268   & 480   & 340   & 370   & 304   & 291   & 369 \\
    352\_3\_5 & \Ffmax & 1213  & 8333  & 1658  & 1133  & 863   & 697   & 781   & 355\_3\_5 & 24592 & 624   & 753   & 457   & 744   & 448   & 388   & 428 \\
    352\_4\_2 & 48585 & 603   & 818   & 679   & 775   & 668   & 460   & 528   & 355\_4\_2 & 363   & 214   & 268   & 226   & 261   & 261   & 203   & 221 \\
    352\_4\_3 & 79649 & 867   & 628   & 720   & 876   & 590   & 470   & 511   & 355\_4\_3 & 714   & 463   & 360   & 369   & 343   & 383   & 260   & 314 \\
    352\_4\_4 & \Ffmax & \Fincr & 4570  & 1046  & 1200  & 858   & 708   & 804   & 355\_4\_4 & 32137 & 404   & 700   & 411   & 532   & 562   & 367   & 451 \\
    
    \bottomrule
    \end{tabular}%

\caption{Results for each system of the sequences generated in the curve section of the train track with velocity ${v}=16\ m/s$.}\label{tab:v16_cur}
\end{sidewaystable}

\end{document}